\documentclass{scrartcl}

\usepackage{amsfonts}
\usepackage{amsmath}
\usepackage{algorithm}
\usepackage{algorithmic}
\usepackage{hyperref}
\usepackage{graphicx}

\newcommand{\bbbr}{\mathbb{R}}
\newcommand{\Idx}{\mathcal{I}}
\newcommand{\Jdx}{\mathcal{J}}
\newcommand{\ctI}{\mathcal{T}_\Idx}
\newcommand{\ctJ}{\mathcal{T}_\Jdx}
\newcommand{\ctIalpha}{\mathcal{T}_{\Idx_\alpha}}
\newcommand{\ctIbeta}{\mathcal{T}_{\Idx_\beta}}
\newcommand{\ctJalpha}{\mathcal{T}_{\Jdx_\alpha}}
\newcommand{\ctJbeta}{\mathcal{T}_{\Jdx_\beta}}
\newcommand{\ctalphabeta}{\mathcal{T}_{\Idx_\alpha\times\Jdx_\beta}}
\newcommand{\ctbetaalpha}{\mathcal{T}_{\Idx_\beta\times\Jdx_\alpha}}
\newcommand{\ctIJ}{\mathcal{T}_{\Idx\times\Jdx}}
\newcommand{\lfI}{\mathcal{L}_\Idx}
\newcommand{\lfJ}{\mathcal{L}_\Jdx}
\newcommand{\lfIJ}{\mathcal{L}_{\Idx\times\Jdx}}
\newcommand{\lfaIJ}{\mathcal{L}_{\Idx\times\Jdx}^+}
\newcommand{\lfiIJ}{\mathcal{L}_{\Idx\times\Jdx}^-}
\newcommand{\nI}{n_\Idx}
\newcommand{\nJ}{n_\Jdx}
\newcommand{\supp}{\mathop{\operatorname{supp}}\nolimits}
\newcommand{\treeroot}{\mathop{\operatorname{root}}\nolimits}
\newcommand{\chil}{\mathop{\operatorname{chil}}\nolimits}
\newcommand{\diam}{\mathop{\operatorname{diam}}\nolimits}
\newcommand{\dist}{\mathop{\operatorname{dist}}\nolimits}
\newcommand{\shareholders}{\mathop{\operatorname{shareholders}}\nolimits}
\newcommand{\manager}{\mathop{\operatorname{manager}}\nolimits}

\newtheorem{definition}{Definition}
\newtheorem{remark}[definition]{Remark}

\begin{document}

\title{\texorpdfstring{Distributed $\mathcal{H}^2$-Matrices for
    Boundary Element Methods}
  {Distributed H2-Matrices for Boundary Element Methods}}

\author{Steffen B\"orm}

\maketitle

\begin{abstract}
  \noindent
  Standard discretization techniques for boundary integral equations,
  e.g., the Galerkin boundary element method, lead to large densely
  populated matrices that require fast and efficient compression techniques
  like the fast multipole method or hierarchical matrices.
  If the underlying mesh is very large, running the corresponding
  algorithms on a distributed computer is attractive, e.g., since
  distributed computers frequently are cost-effective and
  offer a high accumulated memory bandwidth.

  Compared to the closely related particle methods, for which
  distributed algorithms are well-established, the Galerkin
  discretization poses a challenge, since the supports of the
  basis functions influence the block structure of the matrix
  and therefore the flow of data in the corresponding algorithms.
  This article introduces distributed $\mathcal{H}^2$-matrices, a
  class of hierarchical matrices that is closely related to fast
  multipole methods and particularly well-suited for distributed
  computing.
  While earlier efforts required the global tree structure of the
  $\mathcal{H}^2$-matrix to be stored in every node of the distributed
  system, the new approach needs only local multilevel information
  that can be obtained via a simple distributed algorithm, allowing
  us to scale to significantly larger systems.
  Experiments show that this approach can handle very large meshes
  with more than $130$ million triangles efficiently.
\end{abstract}

% ----------------------------------------
% Introduction
% ----------------------------------------
\section{Introduction}

Let us consider the boundary integral equation
\begin{align}\label{eq:integral}
  \int_{\partial\Omega} g(x,y)\, u(y) \,dy
  &= f(x) &
  &\text{ for all } x\in\partial\Omega,
\end{align}
where $\Omega\subseteq\bbbr^3$ is a Lipschitz domain, $g$ is a
kernel function, e.g., the Laplace kernel
\begin{align*}
  g(x,y) &= \frac{1}{4\pi\|x-y\|} &
  &\text{ for all } x,y\in\bbbr^3,\ x\neq y,
\end{align*}
or the Helmholtz kernel
\begin{align*}
  g(x,y) &= \frac{e^{i \kappa \|x-y\|}}{4\pi \|x-y\|} &
  &\text{ for all } x,y\in\bbbr^3,\ x\neq y,
\end{align*}
for the wave number $\kappa\in\bbbr_{\geq 0}$, $f$ is a given
function, and $u$ is the solution we want to compute or at least
approximate.
Integral equations of this kind appear, e.g., when dealing with
exterior domain problems in electrostatics or acoustics.

Galerkin's method offers an elegant approach to discretizing
the integral equation:
we choose a family $(\varphi_i)_{i\in\Idx}$ of \emph{test functions}
and a family $(\psi_j)_{j\in\Jdx}$ of \emph{trial functions} and
look for an approximate solution of the form
\begin{equation}\label{eq:u_h_def}
  u_h := \sum_{j\in\Jdx} z_j \psi_j
\end{equation}
with a coefficient vector $z\in\bbbr^\Jdx$.
Replacing $u$ by $u_h$ in (\ref{eq:integral}), multiplying by
the trial functions $\varphi_i$, and integrating yields
\begin{align}\label{eq:discrete}
  \int_{\partial\Omega} \varphi_i(x)
  \int_{\partial\Omega} g(x,y)\, u_h(y) \,dy \,dx
  &= \int_{\partial\Omega} \varphi_i(x)\, f(x) \,dx &
  &\text{ for all } i\in\Idx,
\end{align}
and with (\ref{eq:u_h_def}) we obtain
\begin{align*}
  \sum_{j\in\Jdx} z_j
     \underbrace{\int_{\partial\Omega} \varphi_i(x) \int_{\partial\Omega} g(x,y)
     \, \psi_j(y) \,dy \,dx}_{=:g_{ij}}
  &= \underbrace{\int_{\partial\Omega} \varphi_i(x)\, f(x) \,dx}_{=:b_i} &
  &\text{ for all } i\in\Idx.
\end{align*}
If we collect the values $g_{ij}$ in a matrix $G\in\bbbr^{\Idx\times\Jdx}$
and the values $b_i$ in a vector $b\in\bbbr^\Idx$, we arrive at the
linear system $Gz=b$.
If we can solve this system, we can use (\ref{eq:u_h_def}) to obtain
the approximate solution $u_h$ of the integral equation (\ref{eq:integral}).

Solving the system $Gz=b$ poses a significant challenge, since
\begin{itemize}
  \item the matrix $G$ is generally not sparse,
  \item its condition number tends to be large, and
  \item computing some of its entries involves evaluating singular
    integrals.
\end{itemize}
The third issue can be addressed by suitable quadrature techniques
\cite{SA96,ERSA98,SASC11,BOHA03a}.
For the second issue, Krylov methods \cite{HA94,FI96,SA03} can be
combined with suitable preconditioners \cite{PEST92,STWE98,LAPURE03,
FAMEPR13b}.

The first issue, i.e., handling large non-sparse matrices efficiently,
is the topic of this article.
If the geometry is locally smooth, wavelet methods
\cite{DAPRSC94b,DASC99,TA02,DAHASC02,HASC06} can be used to obtain
a sparse approximation of the matrix $G$ by taking advantage of
local smoothness properties of the kernel function $g$.
In more general settings, it is usually more attractive to either
approximate the kernel function $g$ directly, e.g., by fast multipole
methods \cite{RO85,GRRO87,AN92,GIRO02}, interpolation
\cite{GI01,BOGR02,BOHA02a,BOLOME02}, quadrature \cite{BOGO12,BOCH14},
or Taylor expansion \cite{HANO89,SA00}, or the use matrix coefficients
to construct local low-rank approximations \cite{TY99,BE00a,BERJ01,BE08}.

Both the approximation of the kernel function and the approximation
of matrix entries leads to a decomposition of the matrix $G$ into
submatrices that are either small enough to be stored directly or
far enough from the singularity at $x=y$ to be approximated by a
low-rank matrix that can be stored efficiently in factorized form.
This leads to what is known as a \emph{hierarchical matrix}
\cite{HA99,HA15}, sometimes even to the sub-class of
$\mathcal{H}^2$-matrices \cite{HAKHSA00,BOHA02,BO10}.

Although hierarchical matrices, and particularly $\mathcal{H}^2$-matrices,
tend to be very efficient, handling very large systems still requires
a large amount of memory and leads to long run-times.
This article addresses this problem by using a distributed computing
approach where multiple computers are connected by a communication
network allowing them to share data and computational work.

In the context of particle systems, several very efficient parallel
algorithms already exist \cite{GRGR90,WASA92,BIYI12,OHYO14}, even for the
particularly challenging Helmholtz problem \cite{LIHU20}.
Generalizing these techniques to boundary element matrices is not
straightforward, since we have to take the shapes of the supports of
the basis functions into account in order to guarantee reliable
convergence.

$\mathcal{H}^2$-matrices can easily deal with general (localized) basis
functions due to their more flexible tree structures \cite{BEBO08}, but
distributing these tree structures among multiple nodes in a distributed
system poses a challenge.
This article presents a technique that can be used to distribute
the trees among all nodes, providing each node with just enough
information to carry out the required operations.

Since the geometry and the choice of basis functions heavily
influence the matrix structure and since we want each node to be
responsible for only part of the domain, the nodes have to
communicate in order to determine what data they need to send to and
receive from each other.
This article presents an efficient algorithm that uses only 
local information to construct \emph{send trees} and \emph{receive trees}
(closely related to the \emph{locally essential trees} used in
parallel multipole methods \cite{WASA92}) that describe the
communication patterns required to carry out operations like
the matrix setup, matrix compression, or matrix-vector multiplication.

Once the send and receive trees are constructed, performing these
operations is fairly straightforward.
This article contains examples of how fundamental operations can
be implemented in a distributed system and provides experimental
results showing weak scaling up to surfaces consisting of more
than $134$ million triangles.

The article is organized as follows:
Section~\ref{se:h2matrices} gives a brief introduction into
the structure and properties of $\mathcal{H}^2$-matrices,
Section~\ref{se:distributed} describes a simple approach to
handling $\mathcal{H}^2$-matrices on a distributed system,
and Section~\ref{se:shared} outlines how this approach can be
improved to obtain better theoretical scaling for very large numbers of
nodes.

% ----------------------------------------
% H^2-matrices
% ----------------------------------------
\section{\texorpdfstring{$\mathcal{H}^2$-matrices}
                        {H2-matrices}}
\label{se:h2matrices}

Before we can discuss how to implement $\mathcal{H}^2$-matrices
on a distributed system, we have to briefly recall how
$\mathcal{H}^2$-matrices work.
In the context of integral equations, interpolation offers a simple
approach:
we fix axis-parallel boxes $\tau,\sigma\subseteq\bbbr^3$ and
denote interpolation points in these boxes by
$(\xi_{\tau,\nu})_{\nu=1}^k$ and $(\xi_{\sigma,\mu})_{\mu=1}^k$.
With the corresponding Lagrange polynomials
$(\ell_{\tau,\nu})_{\nu=1}^k$ and $(\ell_{\sigma,\mu})_{\mu=1}^k$,
the interpolating polynomial of the kernel function $g$ is
given by
\begin{align}\label{eq:interpolation}
  \tilde g_{\tau\sigma}(x,y)
  &= \sum_{\nu=1}^k \sum_{\mu=1}^k g(\xi_{\tau,\nu}, \xi_{\sigma,\mu})\,
       \ell_{\tau,\nu}(x)\, \ell_{\sigma,\mu}(y) &
  &\text{ for all } x\in\tau,\ y\in\sigma.
\end{align}
If $\tau$ and $\sigma$ are well-separated, $g|_{\tau\times\sigma}$ is
smooth and we have
\begin{align*}
  g(x,y) &\approx \tilde g_{\tau\sigma}(x,y) &
  &\text{ for all } x\in\tau,\ y\in\sigma.
\end{align*}
We can take advantage of this approximation of the kernel function
to obtain an approximation of the matrix: we denote by
\begin{align}\label{eq:cluster_indices}
  \hat\tau &:= \{ i\in\Idx\ :\ \supp\varphi_i \subseteq \tau \}, &
  \hat\sigma &:= \{ j\in\Jdx\ :\ \supp\psi_j \subseteq \sigma \}
\end{align}
the indices of basis functions supported in $\tau$ and $\sigma$.
Since $\tilde g_{\tau\sigma}$ approximates $g$ in $\tau\times\sigma$,
we find
\begin{align*}
  g_{ij} &= \int_{\partial\Omega} \varphi_i(x)
    \int_{\partial\Omega} g(x,y)\, \psi_j(y) \,dy\,dx
   \approx \int_{\partial\Omega} \varphi_i(x)
    \int_{\partial\Omega} \tilde g_{\tau\sigma}(x,y)\, \psi_j(y) \,dy\,dx\\
  &= \sum_{\nu=1}^k \sum_{\mu=1}^k
     \underbrace{\int_{\partial\Omega} \varphi_i(x)\, \ell_{\tau,\nu}(x) \,dx
                }_{=:v_{\tau,i\nu}}
     \underbrace{g(\xi_{\tau,\nu},\xi_{\sigma,\mu})}_{=:s_{\tau\sigma,\nu\mu}}
     \underbrace{\int_{\partial\Omega} \psi_j(y)\, \ell_{\sigma,\mu}(y) \,dy
                }_{=:w_{\sigma,j\mu}}
  \quad\text{ for all } i\in\hat\tau,\ j\in\hat\sigma.
\end{align*}
Introducing the matrices $V_\tau\in\bbbr^{\hat\tau\times k}$,
$W_\sigma\in\bbbr^{\hat\sigma\times k}$ (with the shorthand notation
$\bbbr^{\hat\tau\times k}$ for matrices with $k$ columns and row
indices in $\hat\tau$), and $S_{\tau\sigma}\in\bbbr^{k\times k}$ with
\begin{subequations}
\begin{align}
  v_{\tau,i\nu}
  &= \int_{\partial\Omega} \varphi_i(x)\, \ell_{\tau,\nu}(x) \,dx &
  &\text{ for all } i\in\hat\tau,\ \nu\in[1:k],\label{eq:leafV}\\
  w_{\sigma,j\mu}
  &= \int_{\partial\Omega} \psi_j(y)\, \ell_{\sigma,\mu}(y) \,dy &
  &\text{ for all } j\in\hat\sigma,\ \mu\in[1:k],\label{eq:leafW}\\
  s_{\tau\sigma,\nu\mu}
  &= g(\xi_{\tau,\nu}, \xi_{\sigma,\mu}) &
  &\text{ for all } \nu,\mu\in[1:k],
  \label{eq:coupling}
\end{align}
\end{subequations}
we obtain an approximation of the matrix block
\begin{equation}\label{eq:vsw}
  G|_{\hat\tau\times\hat\sigma}
  \approx V_\tau S_{\tau\sigma} W_\sigma^T.
\end{equation}
Unfortunately, we cannot apply this procedure globally, since
$\tau$ and $\sigma$ have to be well-separated.
The solution is to cover $\partial\Omega\times\partial\Omega$ with
multiple well-separated subdomains and a few small remainders that
can be treated directly.
In order to make the construction of this covering as efficient
as possible, we impose a hierarchical structure.

%
% Definition: Cluster tree
%
\begin{definition}[Cluster tree]
\label{de:cluster_tree}
Let $\ctI$ be a tree with root $\rho=\treeroot(\ctI)$, and assign a set
$\hat\tau\subseteq\Idx$ to each of its vertices $\tau\in\ctI$.
We denote the children of a vertex $\tau\in\ctI$ by
$\chil(\tau)\subseteq\ctI$.

$\ctI$ is called a \emph{cluster tree} for the index set $\Idx$
(and the basis functions $(\varphi_i)_{i\in\Idx}$) if
\begin{itemize}
  \item $\hat\varrho = \Idx$,
  \item for all $\tau\in\ctI$, $i\in\hat\tau$, we have
        $\supp\varphi_i\subseteq\tau$,
  \item for all $\tau\in\ctI$ and $\tau_1,\tau_2\in\chil(\tau)$,
        we have $\tau_1\neq\tau_2 \Rightarrow
        \hat\tau_1\cap\hat\tau_2=\emptyset$, and
  \item for all $\tau\in\ctI$ with $\chil(\tau)\neq\emptyset$,
        we have $\hat\tau = \bigcup\{ \hat\tau'\ :\ \tau'\in\chil(\tau) \}$.
\end{itemize}
We denote the leaves of a cluster tree $\ctI$ by $\lfI$.
The vertices of a cluster tree are called \emph{clusters}.
\end{definition}

A simple construction of cluster trees is based on \emph{characteristic
points} $x_i\in\supp\varphi_i$ assigned to represent every basis
function: we recursively split the ``cloud'' of characteristic points,
e.g., by bisection along suitable planes, to obtain a tree structure with
index sets $\hat\tau\subseteq\Idx$ assigned to each vertex.
The clusters $\tau$ can then be constructed by finding minimal axis-parallel
boxes ensuring (\ref{eq:cluster_indices}).

%
% Remark: Cluster geometry
%
\begin{remark}[Cluster geometry]
Our approach requires us to ensure $\supp\varphi_i\subseteq\tau$
for all $i\in\hat\tau$, and this means that for the commonly-used
unstructured surface meshes the clusters $\tau$ will vary in
diameter and aspect ratio, even within the same level of the cluster
tree.

This is a marked difference from classical fast multipole methods
for particle system that typically use a very regular box tree in
order to keep the implementation simple.
We do not have this option, since we are working with entire supports
of basis functions instead of points.

In a distributed system, this implies that nodes cannot simply predict
the clusters in other nodes, since these clusters depend on the local
geometry of the mesh that is not available to all nodes.
\end{remark}

In the following, we assume that cluster trees $\ctI$ for the
basis $(\varphi_i)_{i\in\Idx}$ and $\ctJ$ for the basis
$(\psi_j)_{j\in\Jdx}$ are given.
Our task is now to cover $\partial\Omega\times\partial\Omega$
by pairs of clusters $\tau\times\sigma$ with $\tau\in\ctI$,
$\sigma\in\ctJ$.
We do this by recursively constructing a \emph{block tree}:
the root of the block tree $\ctIJ$ is the pair of the roots
of $\ctI$ and $\ctJ$, i.e., $\treeroot(\ctIJ) = (\treeroot(\ctI),
\treeroot(\ctJ))$.

Given a pair $(\tau,\sigma)\in\ctIJ$, we check whether $\tau$
and $\sigma$ are well-separated, e.g., via the \emph{standard
admissibility condition}
\begin{equation}\label{eq:admissibility}
  \max\{\diam(\tau),\diam(\sigma)\}
  \leq 2 \eta \dist(\tau,\sigma),
\end{equation}
where $\eta\in\bbbr_{>0}$ is a parameter that allows us to balance
the speed of convergence against the storage requirements:
if $\eta$ is small, we have fast convergence, if $\eta$ is large,
we have a small number of blocks.
If $\tau$ and $\sigma$ are well-separated, we make $(\tau,\sigma)$
a leaf of the block tree, i.e., we have $\chil(\tau,\sigma)=\emptyset$.

Otherwise we check whether $\tau$ and $\sigma$ have children.
If one of them does, we let
\begin{equation*}
  \chil(\tau,\sigma) := \begin{cases}
    \chil(\tau)\times\chil(\sigma)
    &\text{ if } \chil(\tau)\neq\emptyset,\ \chil(\sigma)\neq\emptyset,\\
    \{\tau\}\times\chil(\sigma)
    &\text{ if } \chil(\tau)=\emptyset,\ \chil(\sigma)\neq\emptyset,\\
    \chil(\tau)\times\{\sigma\}
    &\text{ if } \chil(\tau)\neq\emptyset,\ \chil(\sigma)=\emptyset,
  \end{cases}
\end{equation*}
and apply our procedure recursively to the children.

If $\chil(\tau)=\emptyset=\chil(\sigma)$, we again make $(\tau,\sigma)$
a leaf of the block tree $\ctIJ$, but this time it is a leaf that does not
allow us to approximate the corresponding matrix block.

We denote the set of leaves of the block tree $\ctIJ$ by $\lfIJ$
and split it into the subset of \emph{admissible leaves} $\lfaIJ$
containing well-separated clusters and the remainder $\lfiIJ$
of \emph{inadmissible leaves}.
Due to the properties of the cluster trees $\ctI$ and $\ctJ$,
we have that $\{ \hat\tau\times\hat\sigma\ :\ (\tau,\sigma)\in\lfIJ \}$
is a partition of $\Idx\times\Jdx$, corresponding to a decomposition
of the matrix $G$ into submatrices $G|_{\hat\tau\times\hat\sigma}$.
For admissible leaves, these submatrices can be approximated in
the factorized form (\ref{eq:vsw}), while inadmissible leaves require
us to store $G|_{\hat\tau\times\hat\sigma}$ directly.

Under standard assumptions, this approach would already reduce
the storage requirements from $\mathcal{O}(\nI\nJ)$ to
$\mathcal{O}((\nI+\nJ) k \log(\nI+\nJ))$, where
$\nI:=|\Idx|$ and $\nJ:=|\Jdx|$ denote the cardinalities of the
index sets and $k$ is the number of interpolation points.

If we use the same interpolation order for all clusters, we can
reduce the complexity even further:
let $\tau\in\ctI$ and $\tau'\in\chil(\tau)$.
Since the interpolation orders for $\tau$ and $\tau'$ are identical,
the Lagrange polynomials of $\tau$ can be expressed by the
Lagrange polynomials of $\tau'$, i.e., we have
\begin{align*}
  \ell_{\tau,\nu}
  &= \sum_{\nu'=1}^k \ell_{\tau,\nu}(\xi_{\tau',\nu'})\,
                    \ell_{\tau',\nu'} &
  &\text{ for all } \nu\in[1:k].
\end{align*}
Introducing the \emph{transfer matrix} $E_{\tau'}\in\bbbr^{k\times k}$
by
\begin{align}\label{eq:transfer}
  e_{\tau',\nu'\nu} &:= \ell_{\tau,\nu}(\xi_{\tau',\nu'}) &
  &\text{ for all } \nu,\nu'\in[1:k],
\end{align}
we obtain
\begin{align*}
  \ell_{\tau,\nu} &= \sum_{\nu'=1}^k e_{\tau',\nu'\nu}\, \ell_{\tau',\nu'} &
  &\text{ for all } \nu\in[1:k],
\end{align*}
and discretizing this equation yields
\begin{align*}
  v_{\tau,i\nu}
  &= \int_{\partial\Omega} \varphi_i(x)\, \ell_{\tau,\nu}(x) \,dx\\
  &= \sum_{\nu'=1}^k e_{\tau',\nu'\nu} \int_{\partial\Omega}
        \varphi_i(x)\, \ell_{\tau',\nu'}(x) \,dx
   = \sum_{\nu'=1}^k e_{\tau',\nu'\nu} v_{\tau',i\nu'} &
  &\text{ for all } i\in\hat\tau',\ \nu\in[1:k].
\end{align*}
This equation can be written in the compact form
\begin{equation}\label{eq:V_nested}
  V_\tau|_{\hat\tau'\times[1:k]} = V_{\tau'} E_{\tau'}.
\end{equation}
Since $\{\hat\tau'\ :\ \tau'\in\chil(\tau)\}$ is a partition of
$\hat\tau$, we can use (\ref{eq:V_nested}) to reconstruct the
entire matrix $V_\tau$ using the matrices $V_{\tau'}$ corresponding
to the children $\tau'$ of $\tau$ and the small $k\times k$
transfer matrices.
We can apply this procedure recursively and see that we have
to store $V_\tau$ only for leaves $\tau\in\lfI$, while
the transfer matrices can be used for all other clusters.
Treating $W_\sigma$ in the same way (using the notation $F_\sigma$ for the
transfer matrices) reduces the storage complexity
to $\mathcal{O}(\nI k + \nJ k)$ \cite{BO10}, i.e., we have the optimal
order of complexity if $k$ is constant.
In practice, we have to use $k\sim\log(\nI+\nJ)^3$ to keep the
accuracy of the matrix approximation consistent with the
discretization error.
Further improvements are possible \cite{SA00,BOLOME02,BOSA03},
but are not the topic of this brief introduction.

The matrix $G$ is now approximated by the \emph{coupling matrices}
$S_{\tau\sigma}\in\bbbr^{k\times k}$ for admissible blocks
$(\tau,\sigma)\in\lfaIJ$, the \emph{nearfield matrices}
$G|_{\hat\tau\times\hat\sigma}$ for inadmissible blocks
$(\tau,\sigma)\in\lfiIJ$, the \emph{leaf matrices}
$V_\tau$ and $W_\sigma$ for $\tau\in\lfI$ and $\sigma\in\lfJ$,
and the \emph{transfer matrices} $E_\tau,F_\sigma\in\bbbr^{k\times k}$
for all $\tau\in\ctI$, $\sigma\in\ctJ$.
In order to perform a matrix-vector multiplication, i.e., to
compute $y \gets y + G x$ with $x\in\bbbr^\Jdx$ and $y\in\bbbr^\Idx$,
we have to perform updates
\begin{align*}
  y|_{\hat\tau} &\gets y|_{\hat\tau}
    + \begin{cases}
        V_\tau S_{\tau\sigma} W_\sigma^T x|_{\hat\sigma}
        &\text{ if } (\tau,\sigma)\in\lfaIJ,\\
        G|_{\hat\tau\times\hat\sigma} &\text{ otherwise}
      \end{cases} &
  &\text{ for all } (\tau,\sigma)\in\lfIJ.
\end{align*}
We split the computation into three phases:
\begin{enumerate}
  \item \emph{Forward transformation:}
    we compute $\hat x_\sigma := W_\sigma^T x|_{\hat\sigma}$
    for all $\sigma\in\ctJ$.
  \item \emph{Interaction:}
    for all admissible blocks $(\tau,\sigma)\in\lfaIJ$, we add
    $S_{\tau\sigma} \hat x_\sigma$ to $\hat y_\tau$.

    For all inadmissible blocks $(\tau,\sigma)\in\lfiIJ$, we add
    $G|_{\hat\tau\times\hat\sigma} x|_{\hat\sigma}$ to
    $y|_{\hat\tau}$.
  \item \emph{Backward transformation:}
    we add $V_\tau \hat y_\tau$ to $y|_{\hat\tau}$ for all $\tau\in\ctI$.
\end{enumerate}
In order to get the correct result, we have to set
$\hat y_\tau := 0\in\bbbr^k$ for all $\tau\in\ctI$ before starting
the interaction phase of the algorithm.

Using the transfer matrices allows us to obtain a very elegant
algorithm for the forward transformation:
if $\sigma\in\ctJ$ is a leaf, we compute $\hat x_\sigma = W_\sigma^T
x|_{\hat\sigma}$ directly.
Otherwise, we ensure that all children $\sigma'\in\chil(\sigma)$ are treated
\emph{before} $\sigma$, so that $\hat x_{\sigma'}$ is available, and then
use (\ref{eq:V_nested}) to get
\begin{equation*}
  \hat x_\sigma
  = W_\sigma^T x|_{\hat\sigma}
  = \sum_{\sigma'\in\chil(\sigma)}
      F_{\sigma'}^T W_{\sigma'}^T x|_{\hat\sigma'}
  = \sum_{\sigma'\in\chil(\sigma)}
      F_{\sigma'}^T \hat x_{\sigma'}
\end{equation*}
to compute $\hat x_\sigma$ efficiently, cf. Algorithm~\ref{al:forward}.

%
% Algorithm: Forward transformation
%
\begin{algorithm}
  \caption{Forward transformation}
  \label{al:forward}
  \begin{algorithmic}[1]
    \STATE \textbf{procedure} forward($\sigma$, $x$, \textbf{var} $\hat x$);
    \IF{$\chil(\sigma)=\emptyset$}
      \STATE $\hat x_\sigma \gets W_\sigma^T x|_{\hat\sigma}$
    \ELSE
      \STATE $\hat x_\sigma \gets 0$;
      \FOR{$\sigma'\in\chil(\sigma)$}
        \STATE forward($\sigma'$, $x$, $\hat x$);
        \STATE $\hat x_\sigma \gets \hat x_\sigma + F_{\sigma'}^T \hat x_{\sigma'}$
      \ENDFOR
    \ENDIF
  \end{algorithmic}
\end{algorithm}

We can use a similar approach for the backward transformation:
if $\tau\in\ctI$ is a leaf, we add $V_\tau \hat y_\tau$ directly
to $y|_{\hat\tau}$.
Otherwise, we notice that due to (\ref{eq:V_nested}) the update
$y|_{\hat\tau} \gets y|_{\hat\tau} + V_\tau \hat y_\tau$ is
equivalent with
\begin{align*}
  y|_{\hat\tau'} &\gets y|_{\hat\tau'} + V_{\tau'} E_{\tau'} \hat y_\tau &
  &\text{ for all } \tau'\in\chil(\tau).
\end{align*}
This means that instead of adding $V_\tau \hat y_\tau$ directly to
$y|_{\hat\tau}$, we can simply add $E_{\tau'} \hat y_\tau$ to
$\hat y_{\tau'}$ instead, as long as we ensure that the children
$\tau'\in\chil(\tau)$ are treated \emph{after} $\tau$.
This can be ensured by a simple recursion, cf. Algorithm~\ref{al:backward}.

%
% Algorithm: Backward transformation
%
\begin{algorithm}
  \caption{Backward transformation}
  \label{al:backward}
  \begin{algorithmic}[1]
    \STATE \textbf{procedure} backward($\tau$, \textbf{var} $\hat y$, $y$);
    \IF{$\chil(\tau)=\emptyset$}
      \STATE $y|_{\hat\tau}\gets y|_{\hat\tau} + V_\tau \hat y_\tau$
    \ELSE
      \FOR{$\tau'\in\chil(\tau)$}
        \STATE $\hat y_{\tau'} \gets \hat y_{\tau'} + E_{\tau'} \hat y_\tau$;
        \STATE backward($\tau'$, $\hat y$, $y$)
      \ENDFOR
    \ENDIF
  \end{algorithmic}
\end{algorithm}

%
% Algorithm: Interaction phase
%
\begin{algorithm}
  \caption{Interaction phase}
  \label{al:interaction}
  \begin{algorithmic}[1]
    \STATE \textbf{procedure} interaction($(\tau,\sigma)$, $\hat x$, $x$,
       \textbf{var} $\hat y$, $y$);
    \IF{$\chil(\tau,\sigma)\neq\emptyset$}
       \FOR{$(\tau',\sigma')\in\chil(\tau,\sigma)$}
          \STATE interaction($(\tau',\sigma')$, $\hat x$, $x$, $\hat y$, $y$)
       \ENDFOR
    \ELSIF{$(\tau,\sigma)$ is admissible}
       \STATE $\hat y_\tau \gets \hat y_\tau + S_{\tau\sigma} \hat x_\sigma$
    \ELSE
       \STATE $y|_\tau \gets y|_\tau + G|_{\hat\tau\times\hat\sigma}
                  x|_\sigma$
    \ENDIF
  \end{algorithmic}
\end{algorithm}

The interaction phase can also be handled by recursion,
cf. Algorithm~\ref{al:interaction}, and combining these algorithms,
the matrix-vector multiplication can be performed in
$\mathcal{O}(\nI k + \nJ k)$ operations \cite{BO10}.

% ----------------------------------------
% Distributed H2-matrices
% ----------------------------------------
\section{\texorpdfstring{Distributed $\mathcal{H}^2$-matrices}
                        {Distributed H2-matrices}}
\label{se:distributed}

In order to be able to treat very large matrices, we aim to
use a distributed system of $p$ processing nodes connected
by a fast network.
The indices $\Idx$ of the test functions and $\Jdx$ of the trial
functions are split among the nodes, preferably according to suitable
subdomains, and we denote the subsets of indices for a node
$\alpha\in[1:p]$ by $\Idx_\alpha\subseteq\Idx$ and $\Jdx_\alpha\subseteq\Jdx$.
These subsets have to form partitions, i.e.,
\begin{gather*}
  \Idx_\alpha \cap \Idx_\beta = \emptyset,\qquad
  \Jdx_\alpha \cap \Jdx_\beta = \emptyset
  \qquad\text{ for all } \alpha,\beta\in[1:p] \text{ with } \alpha\neq\beta,\\
  \Idx = \bigcup_{\alpha=1}^p \Idx_\alpha,\qquad
  \Jdx = \bigcup_{\alpha=1}^p \Jdx_\alpha.
\end{gather*}
In order to use each nodes' storage efficiently, the supports
\begin{align*}
  \Gamma_\alpha &:= \bigcup\{ \supp\varphi_i,\ \supp\psi_j
                      \ :\ i\in\Idx_\alpha,\ j\in\Jdx_\alpha \} &
  &\text{ for all } \alpha\in[1:p]
\end{align*}
of all basis functions assigned to a node $\alpha\in[1:p]$ should
intersect as little as possible with the supports of other nodes.
The matrix $G$ can now be written in block form as
\begin{equation*}
  G = \begin{pmatrix}
        G|_{\Idx_1\times\Jdx_1} & \cdots & G|_{\Idx_1\times\Jdx_p}\\
        \vdots & \ddots & \vdots\\
        G|_{\Idx_p\times\Jdx_1} & \cdots & G|_{\Idx_p\times\Jdx_p}
      \end{pmatrix},
\end{equation*}
and we aim to approximate all blocks by $\mathcal{H}^2$-matrices.
To this end, we can apply standard algorithms \cite[Section~5.4]{HA15}
to construct local cluster trees $\ctIalpha$ and $\ctJalpha$ for the index
sets $\Idx_\alpha$ and $\Jdx_\alpha$.
Considering that we want to implement the matrix-vector multiplication,
we decide that a node $\alpha\in[1:p]$ will be responsible for
storing one block row, i.e., the matrix blocks
$G|_{\Idx_\alpha\times\Jdx_1},\ldots,G|_{\Idx_\alpha\times\Jdx_p}$.
In order to treat a block $G|_{\Idx_\alpha\times\Jdx_\beta}$, the
responsible node $\alpha$ needs information stored in node $\beta$,
i.e., clusters in $\ctJbeta$ and geometry information.
In a departure from previous algorithms that stored \emph{all} of this
information in \emph{all} nodes \cite{BEBO08}, we aim to store
the minimum required by our algorithms.
With admissibility conditions like (\ref{eq:admissibility}) in mind,
it is reasonable to expect that if $\Gamma_\alpha$ and $\Gamma_\beta$
are sufficiently far from each other, $G|_{\Idx_\alpha\times\Jdx_\beta}$
can be approximated by very few blocks, possibly even just one,
and therefore we do not have to copy the entire cluster tree $\ctJbeta$
to node $\alpha$.

% ----------------------------------------
% Send and receive trees
% ----------------------------------------
\subsection{Send and receive trees}
A first attempt at an algorithm could look as follows: every node
$\alpha\in[1:p]$ broadcasts the roots of $\ctIalpha$ and $\ctJalpha$
to all other nodes.
Now it checks if $\treeroot(\ctIalpha)$ and $\treeroot(\ctJbeta)$
satisfy the admissibility condition for all $\beta\in[1:p]$.
If the condition holds, the entire block $G|_{\Idx_\alpha\times\Jdx_\beta}$
can be approximated in the form (\ref{eq:vsw}).
Otherwise, the block needs to be split, i.e., node $\alpha$ needs
to get the children of $\treeroot(\ctJbeta)$ from node $\beta$.
If we want to use standard MPI point-to-point communication functions,
this poses a problem: node $\alpha$ knows that it wants to receive data
on the children from node $\beta$, but $\beta$ does not know that it is
supposed to send this data to node $\alpha$.

A simple solution to this problem is to have node $\alpha$ not
only construct block trees $\ctalphabeta$ for the index sets
$\Idx_\alpha\times\Jdx_\beta$, but also block trees $\ctbetaalpha$
for the index sets $\Idx_\beta\times\Jdx_\alpha$ telling us what
information has to be sent from node $\alpha$ to node $\beta$.
If we construct these block trees level-by-level, every node
knows not only what information it needs to receive from other
nodes, but also what information it needs to send.

To keep track of the communication requirements, we introduce
subtrees of cluster trees.
Let $\alpha,\beta\in[1:p]$.
A subtree $\mathcal{T}_\text{col}^{(\alpha,\beta)}\subseteq\ctJbeta$ is called a
\emph{column transmission tree} for the destination $\alpha$ and the
origin $\beta$ if it is the minimal subtree satisfying
\begin{align*}
  \sigma\in\mathcal{T}_\text{col}^{(\alpha,\beta)} &\iff
  \exists \tau\in\ctIalpha\ :\ (\tau,\sigma)\in\ctalphabeta &
  &\text{ for all } \sigma\in\ctJbeta,
\end{align*}
i.e., if it contains exactly the column clusters of $\ctJbeta$
required to construct the block tree $\ctalphabeta$.

A subtree $\mathcal{T}_\text{row}^{(\alpha,\beta)}\subseteq\ctIbeta$ is called a
\emph{row transmission tree} for the destination $\alpha$ and the
origin $\beta$ if it is the minimal subtree satisfying
\begin{align*}
  \tau\in\mathcal{T}_\text{row}^{(\alpha,\beta)} &\iff
  \exists \sigma\in\ctJalpha\ :\ (\tau,\sigma)\in\ctbetaalpha &
  &\text{ for all } \tau\in\ctIbeta,
\end{align*}
i.e., if it contains exactly the row clusters of $\ctIbeta$
required to construct the block tree $\ctbetaalpha$.

At the destination $\alpha$, we call the row and column transmission
trees \emph{receive trees}, while we call them \emph{send trees}
at the origin $\beta$.
We have to ensure that the receive trees constructed at node $\alpha$
exactly match the send trees constructed at the node $\beta$ so that
all nodes agree on what data has to be sent and received.

%
% Figure: Block and transmission trees
%
\begin{figure}
  \begin{center}
    \includegraphics[width=0.9\textwidth]{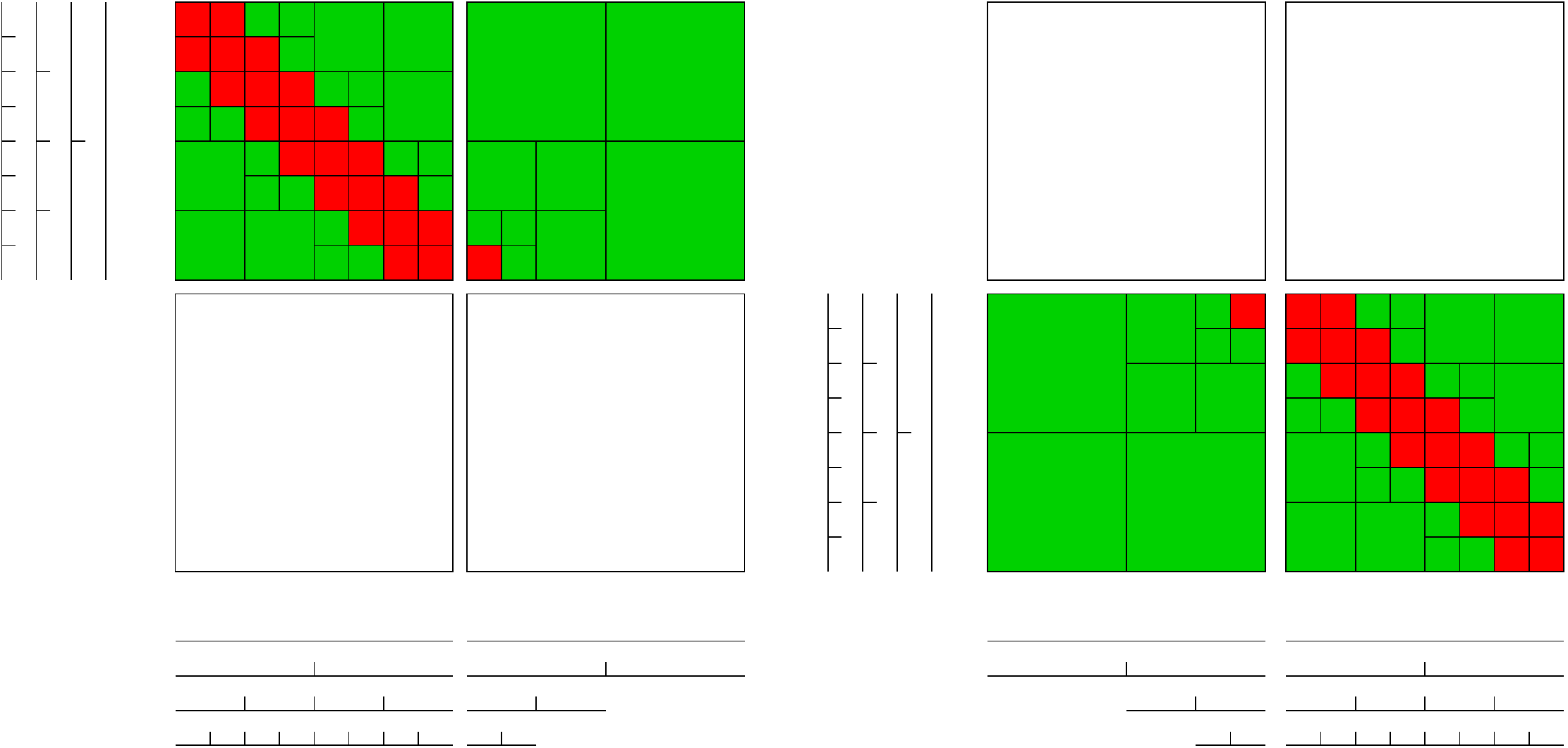}
  \end{center}
  \caption{Local block trees and transmission trees for two nodes.
     Left: Cluster tree $\mathcal{T}_{\Idx_1}$, block trees $\mathcal{T}_{\Idx_1\times\Idx_1}$ and
           $\mathcal{T}_{\Idx_1\times\Idx_2}$ with the column receive
           trees $\mathcal{R}_{\text{col},1}^{(1)}$ and
           $\mathcal{R}_{\text{col},2}^{(1)}$.
     Right: Cluster tree $\mathcal{T}_{\Idx_2}$, block trees $\mathcal{T}_{\Idx_2\times\Idx_1}$ and
           $\mathcal{T}_{\Idx_2\times\Idx_2}$ with the column receive
           trees $\mathcal{R}_{\text{col},1}^{(2)}$ and
           $\mathcal{R}_{\text{col},2}^{(2)}$.}
  \label{fi:sendrecv}
\end{figure}

Figure~\ref{fi:sendrecv} illustrates the connections between the
trees for a simple one-dimensional model problem:
the clusters are intervals, the admissibility condition
(\ref{eq:admissibility}) leads to small blocks near the diagonal
and large blocks far from the diagonal.
Row and column transmission trees are identical, since $\ctI$ and $\ctJ$
are identical.
The transmission trees $\mathcal{T}^{(1,2)}_\text{row}$ and
$\mathcal{T}^{(2,1)}_\text{row}$ contain no more than two clusters per level,
since only these two are needed for the off-diagonal block trees
$\mathcal{T}_{\Idx_1\times\Jdx_2}$ and $\mathcal{T}_{\Idx_2\times\Jdx_1}$.

%
% Algorithm: Construction of block, send, and receive trees
%
\begin{algorithm}
  \caption{Construction of local block trees, send and receive trees}
  \label{al:build_block}
  \begin{algorithmic}[1]
    \STATE \textbf{procedure} build\_block(\textbf{var}
               $\mathcal{T}_{\Idx_\alpha\times\Jdx}$,
               $\mathcal{T}_{\Idx\times\Jdx_\alpha}$,
               $\mathcal{R}^{(\alpha)}_\text{row}$,
               $\mathcal{R}^{(\alpha)}_\text{col}$,
               $\mathcal{S}^{(\alpha)}_\text{row}$,
               $\mathcal{S}^{(\alpha)}_\text{col}$);
    \FOR{$\beta\in[1:p]$}
      \STATE bcast\_cluster($\beta$, $\treeroot(\ctIalpha)$,
                            $\rho_\beta$);\quad
       $\treeroot(\mathcal{R}^{(\alpha)}_{\text{row},\beta}) \gets \rho_\beta$;\quad
       $\treeroot(\mathcal{S}^{(\alpha)}_{\text{row},\beta}) \gets \treeroot(\ctIalpha)$;
      \STATE bcast\_cluster($\beta$, $\treeroot(\ctJalpha)$,
                            $\pi_\beta$);\quad
       $\treeroot(\mathcal{R}^{(\alpha)}_{\text{col},\beta}) \gets \pi_\beta$;\quad
       $\treeroot(\mathcal{S}^{(\alpha)}_{\text{col},\beta}) \gets \treeroot(\ctJalpha)$;
      \STATE $\treeroot(\ctalphabeta) \gets (\treeroot(\ctIalpha),\pi_\beta)$;\quad
             $\mathcal{A}_{\text{row},\beta} \gets \{ (\treeroot(\ctIalpha),\pi_\beta) \}$;
      \STATE $\treeroot(\ctbetaalpha) \gets (\rho_\beta,\treeroot(\ctJalpha))$;\quad
             $\mathcal{A}_{\text{col},\beta} \gets \{ (\rho_\beta,\treeroot(\ctJalpha)) \}$
    \ENDFOR;
    \REPEAT
      \FOR{$\beta\in[1:p]$}
        \STATE $\mathcal{G}_{\text{col},\beta} \gets \emptyset$;\quad
               $\mathcal{P}_{\text{row},\beta} \gets \emptyset$;\quad
               $\mathcal{A}_{\text{row},\beta}^\text{old} \gets
                  \mathcal{A}_{\text{row},\beta}$;\quad
               $\mathcal{A}_{\text{row},\beta} \gets \emptyset$;
        \FOR{$(\tau,\sigma)\in\mathcal{A}_{\text{row},\beta}^\text{old}$}
          \IF{$(\tau,\sigma)$ not admissible and $\chil(\sigma)\neq\emptyset$}
            \STATE $\mathcal{G}_{\text{col},\beta} \gets
                    \mathcal{G}_{\text{col},\beta} \cup \{ \sigma \}$;\quad
                   $\mathcal{P}_{\text{row},\beta} \gets
                    \mathcal{P}_{\text{row},\beta} \cup \{ \tau \}$
          \ENDIF
        \ENDFOR;
        \STATE $\mathcal{G}_{\text{row},\beta} \gets \emptyset$;\quad
               $\mathcal{P}_{\text{col},\beta} \gets \emptyset$;\quad
               $\mathcal{A}_{\text{col},\beta}^\text{old} \gets
                  \mathcal{A}_{\text{col},\beta}$;\quad
               $\mathcal{A}_{\text{col},\beta} \gets \emptyset$;
        \FOR{$(\tau,\sigma)\in\mathcal{A}_{\text{col},\beta}^\text{old}$}
          \IF{$(\tau,\sigma)$ not admissible and $\chil(\tau)\neq\emptyset$}
            \STATE $\mathcal{G}_{\text{row},\beta} \gets
                    \mathcal{G}_{\text{row},\beta} \cup \{ \tau \}$;\quad
                   $\mathcal{P}_{\text{col},\beta} \gets
                    \mathcal{P}_{\text{col},\beta} \cup \{ \sigma \}$
          \ENDIF
        \ENDFOR
      \ENDFOR;
      \STATE Send children of clusters in $\mathcal{P}_{\text{row},\beta}$
             and $\mathcal{P}_{\text{col},\beta}$ to $\beta$,
             copy to $\mathcal{S}^{(\alpha)}_{\text{row},\beta}$
             and $\mathcal{S}^{(\alpha)}_{\text{col},\beta}$;
      \STATE Receive children of clusters in $\mathcal{G}_{\text{row},\beta}$
             and $\mathcal{G}_{\text{col},\beta}$ from $\beta$,
             update $\mathcal{R}^{(\alpha)}_{\text{row},\beta}$
             and $\mathcal{R}^{(\alpha)}_{\text{col},\beta}$;
      \FOR{$\beta\in[1:p]$}
        \FOR{$(\tau,\sigma)\in\mathcal{A}_{\text{row},\beta}^\text{old}$}
          \IF{$(\tau,\sigma)$ not admissible and ($\chil(\tau)\neq\emptyset$
              or $\chil(\sigma)\neq\emptyset$)}
            \STATE Split $(\tau,\sigma)$ and
                   add children to $\mathcal{A}_{\text{row},\beta}$
          \ENDIF
        \ENDFOR
        \FOR{$(\tau,\sigma)\in\mathcal{A}_{\text{col},\beta}^\text{old}$}
          \IF{$(\tau,\sigma)$ not admissible and ($\chil(\tau)\neq\emptyset$
              or $\chil(\sigma)\neq\emptyset$)}
            \STATE Split $(\tau,\sigma)$ and
                   add children to $\mathcal{A}_{\text{col},\beta}$
          \ENDIF
        \ENDFOR
      \ENDFOR;
      \STATE busy=$(\exists \beta\in[1:p]\ :\ 
                    \mathcal{A}_{\text{row},\beta}\neq\emptyset \text{ or }
                    \mathcal{A}_{\text{col},\beta}\neq\emptyset)$;
      \STATE reduce\_or\_bool(busy)
    \UNTIL{busy=false}
  \end{algorithmic}
\end{algorithm}

Algorithm~\ref{al:build_block} is a realization of this approach.
If it is run on a node $\alpha\in[1:p]$, it constructs several
trees central to subsequent algorithms:
\begin{itemize}
  \item block trees $\ctalphabeta$ and $\ctbetaalpha$ for all
        $\beta\in[1:p]$,
  \item receive trees $\mathcal{R}^{(\alpha)}_{\text{row},\beta}
        = \mathcal{T}^{(\alpha,\beta)}_\text{row}$ containing
        copies of exactly those clusters in $\ctIbeta$ that are
        required by node $\alpha$ to build $\ctbetaalpha$,
  \item receive trees $\mathcal{R}^{(\alpha)}_{\text{col},\beta}
        = \mathcal{T}^{(\alpha,\beta)}_\text{col}$ containing
        copies of exactly those clusters in $\ctJbeta$ that are
        required by node $\alpha$ to build $\ctalphabeta$,
  \item send trees $\mathcal{S}^{(\alpha)}_{\text{row},\beta}
        = \mathcal{T}^{(\beta,\alpha)}_\text{row}$ containing
        copies of exactly those clusters in $\ctIalpha$ that are
        in the receive tree $\mathcal{R}^{(\beta)}_{\text{row},\alpha}$ in
        node $\beta$, and
  \item send trees $\mathcal{S}^{(\alpha)}_{\text{col},\beta}
        = \mathcal{T}^{(\beta,\alpha)}_\text{col}$ containing
        copies of exactly those clusters in $\ctJalpha$ that are
        in the receive tree $\mathcal{R}^{(\beta)}_{\text{col},\alpha}$ in
        node $\beta$.
\end{itemize}
Clusters are transmitted in two ways during the course of this
algorithm:
in the initialization phase, the root clusters of the local cluster
trees are broadcast to all nodes.
Sending a cluster $\tau$ means sending only the minimal and maximal
coordinates of the axis-parallel box $\tau$ and the number of
its children, but \emph{not} the children themselves.

In the main loop of the algorithm, we decide whether we need the
children of a cluster $\tau$.
If we do, we again transmit only the coordinates of \emph{all} children
$\tau'\in\chil(\tau)$ and the numbers $\#\chil(\tau')$ of their
children, i.e., the grandchildren, but again not the grandchildren
themselves.
Once this data has been received, the children are created and
suitable pointers in the parent cluster are initialized.

Lines~2 to 7 serve to initialize the root level of the block trees
$\ctalphabeta$ for all $\beta\in[1:p]$:
every node $\beta$ broadcasts the roots $\rho_\beta$ and $\pi_\beta$
of $\ctIbeta$ and $\ctJbeta$, and all nodes use them to construct the
roots of $\ctalphabeta$ and $\ctbetaalpha$.
$\rho_\beta$ becomes the root of the row receive tree
$\mathcal{R}^{(\alpha)}_{\text{row},\beta}$,
$\pi_\beta$ becomes the root of the column receive tree
$\mathcal{R}^{(\alpha)}_{\text{col},\beta}$.
The roots of the row and column send trees
$\mathcal{S}^{(\alpha)}_{\text{row},\beta}$ and
$\mathcal{S}^{(\alpha)}_{\text{col},\beta}$ are set to the roots
of the local cluster trees $\ctIalpha$ and $\ctJalpha$, respectively.
Sets $\mathcal{A}_{\text{row},\beta}$ and $\mathcal{A}_{\text{col},\beta}$
are used to keep track of ``active'' blocks in $\ctalphabeta$ and
$\ctbetaalpha$ that may need to be split if they happen to not be
admissible.

After this initialization phase, we enter the main loop.
Lines~11 to 15 check the admissibility of the active blocks in
$\ctalphabeta$.
If a block has to be split, it column cluster is added to the
set $\mathcal{G}_{\text{col},\beta}$ of clusters whose children we
have to get from node $\beta$, while its row cluster is added
to the set $\mathcal{P}_{\text{row},\beta}$ of clusters whose children
we have to send to node $\beta$.

Lines~16 to 21 perform the same operation for the active blocks
in $\ctbetaalpha$.

Once these checks have been completed for all active blocks, we
know exactly which children have to be sent and received, so we
can carry out the required communication operations in the
lines~23 and 24.
This step can be realized, e.g., by MPI's \verb$Isend$ and \verb$Irecv$
calls or, maybe more efficiently, by a collective \verb$Alltoallv$ call.
After the children have been sent, we can update the send trees
$\mathcal{S}^{(\alpha)}_{\text{row},\beta}$ and
$\mathcal{S}^{(\alpha)}_{\text{col},\beta}$ to keep track of this
communication step.
After the children have been received, we can update the receive
trees $\mathcal{R}^{(\alpha)}_{\text{row},\beta}$ and
$\mathcal{R}^{(\alpha)}_{\text{col},\beta}$ so the children can be
used when splitting inadmissible active blocks.

Lines~26 to 30 perform this operation for blocks in $\ctalphabeta$:
inadmissible blocks are split, as long as they have children, and
the newly created blocks are added to the set $\mathcal{A}_{\text{row},\beta}$
of active blocks for the next iteration.

Lines~31 to 35 do the same for blocks in $\ctbetaalpha$.

Since we generally cannot guarantee that all nodes will stop splitting
blocks at the same time, we have to vote:
a node has to keep working as long as it has active blocks, i.e., as
long as the sets $\mathcal{A}_{\text{row},\beta}$ and
$\mathcal{A}_{\text{col},\beta}$ are not empty.
The entire algorithm has to run for as long as at least one node
still has active blocks, so we use a reduction operation with
the Boolean ``$\operatorname{or}$'' operator to find out if we have to
repeat the loop once more.
This voting process takes place in the lines~37 and 38.

% ----------------------------------------
% Distributed matrix setup
% ----------------------------------------
\subsection{Distributed matrix setup}

Once the block, send, and receive trees have been constructed, we
can set up the matrix itself.
If we base the $\mathcal{H}^2$-matrix on interpolation, this is
very simple:
for all clusters $\tau\in\ctIalpha$ and $\sigma\in\ctJalpha$, we can construct
the transfer matrices $E_\tau$ and $F_\sigma$ directly using
(\ref{eq:transfer}), while the leaf matrices $V_\tau$ and $W_\sigma$ can
be constructed for leaf clusters by using (\ref{eq:leafV}) and
(\ref{eq:leafW}).

Setting up the coupling matrices $S_{\tau\sigma}$ for admissible
blocks $(\tau,\sigma)\in\ctalphabeta$, $\beta\in[1:p]$, is also
straightforward, since we have already transmitted the bounding
boxes of $\sigma$ during the construction of the block tree and
therefore can compute the interpolation points and apply
(\ref{eq:coupling}) directly.

The nearfield matrices $G|_{\hat\tau\times\hat\sigma}$ pose a
somewhat greater challenge, since evaluating the matrix entries
\begin{align*}
  g_{ij} &= \int_{\partial\Omega} \varphi_i(x)
           \int_{\partial\Omega} g(x,y)\, \psi_j(y) \,dy\,dx &
  &\text{ for all } i\in\hat\tau,\ j\in\hat\sigma,
\end{align*}
requires geometric information on the supports of $\varphi_i$
and $\psi_j$, and only information on $\varphi_i$ is locally
available in node $\alpha$.

This problem can be easily solved by transmitting all information
needed for the leaves of the column receive tree
$\mathcal{R}^{(\alpha)}_{\text{col},\beta}$, i.e., the indices of
the basis functions of a given leaf cluster
$\sigma\in\mathcal{R}^{(\alpha)}_{\text{col},\beta}$ and the
geometric information required by the nearfield quadrature
rule for these basis functions.

A similar problem occurs if more advanced approximation techniques
are used.
As an example, we consider the GCA-$\mathcal{H}^2$ method
\cite{BOCH14} that chooses pivot elements $\hat\tau_0\subseteq\hat\tau$
for every $\tau\in\ctIalpha$ and $\hat\sigma_0\subseteq\hat\sigma$
for every $\sigma\in\ctJbeta$ and then uses the submatrix
$G|_{\hat\tau_0\times\hat\sigma_0}$ to compute the coupling matrix
$S_{\tau\sigma}$ for $(\tau,\sigma)\in\ctalphabeta$.
In order to use GCA-$\mathcal{H}^2$, we therefore have to transmit
the geometric information required to evaluate
$G|_{\hat\tau_0\times\hat\sigma_0}$ from node $\beta$ to node $\alpha$.
Fortunately the receive tree $\mathcal{R}^{(\alpha)}_{\text{col},\beta}$
is usually quite small for $\beta\neq\alpha$, and fortunately the
sets $\hat\sigma_0\subseteq\hat\sigma$ chosen by the GCA-$\mathcal{H}^2$
algorithm for $\sigma\in\mathcal{R}^{(\alpha)}_{\text{col},\beta}$ are also
small.
The construction of the subsets $\hat\tau_0\subseteq\hat\tau$
and $\hat\sigma_0\subseteq\hat\sigma$ for $\tau\in\ctIalpha$ and
$\sigma\in\ctJalpha$ can be performed without communication between
the nodes.

% ----------------------------------------
% Distributed matrix-vector multiplication
% ----------------------------------------
\subsection{Distributed matrix-vector multiplication}

The $\mathcal{H}^2$-matrix-vector multiplication consists of three
phases:
\begin{enumerate}
  \item The forward transformation, i.e., the computation of
        $\hat x_\sigma = W_\sigma^T x|_{\hat\sigma}$ for all
        $\sigma\in\ctJalpha$, can be carried out completely locally
        for every node $\alpha\in[1:p]$.
  \item The interaction phase, i.e., the updates
        $\hat y_\tau \gets \hat y_\tau + S_{\tau\sigma} \hat x_\sigma$
        for admissible blocks
        $(\tau,\sigma)\in\mathcal{L}^+_{\Idx_\alpha\times\Jdx_\beta}$ and
        $y|_{\hat\tau} \gets y|_{\hat\tau} + G|_{\hat\tau\times\hat\sigma}
        x|_{\hat\sigma}$ for inadmissible blocks
        $(\tau,\sigma)\in\mathcal{L}^-_{\Idx_\alpha\times\Jdx_\beta}$,
        require us to transmit the auxiliary vectors
        $\hat x_\sigma$ for all
        $\sigma\in\mathcal{R}^{(\alpha)}_{\text{col},\beta}$ and
        the vectors
        $x|_{\hat\sigma}$ for all leaves of
        $\sigma\in\mathcal{R}^{(\alpha)}_{\text{col},\beta}$.
  \item The backward transformation, i.e., the updates
        $y|_{\hat\tau} \gets y|_{\hat\tau} + V_\tau \hat y_\tau$
        for all $\tau\in\ctIalpha$, can again be carried out
        completely locally for every node $\alpha\in[1:p]$.
\end{enumerate}
Algorithm~\ref{al:mvm} shows the complete procedure:
we can use Algorithms~\ref{al:forward} and \ref{al:backward}
for the forward and backward transformation and
Algorithm~\ref{al:interaction} for the interaction phase.
The only new part is the communication phase where we sent
the coefficients of clusters in the column send trees
$\mathcal{S}^{(\alpha)}_{\text{col},\beta}$ and receive the
coefficients of clusters in the column receive trees
$\mathcal{R}^{(\alpha)}_{\text{col},\beta}$ for all $\beta\in[1:p]$.

%
% Algorithm: Distributed matrix-vector multiplication
%
\begin{algorithm}
  \caption{Distributed matrix-vector multiplication}
  \label{al:mvm}
  \begin{algorithmic}[1]

    \STATE \textbf{procedure} mvm(
               $\mathcal{T}_{\Idx_\alpha\times\Jdx}$,
               $\mathcal{R}^{(\alpha)}_\text{col}$,
               $\mathcal{S}^{(\alpha)}_\text{col}$,
               $x$, \textbf{var} y);
    \STATE forward($\treeroot(\ctJalpha)$, $x$, $\hat x$);
    \FOR{$\beta\in[1:p]$}
       \STATE Send $\hat x_\sigma$ to node $\beta$ for all
              $\sigma\in\mathcal{S}^{(\alpha)}_{\text{col},\beta}$;
       \STATE Send $x|_{\hat\sigma}$ to node $\beta$ for all
              $\sigma\in\mathcal{S}^{(\alpha)}_{\text{col},\beta}$
              that are leaves in $\ctJalpha$;
       \STATE Receive $\hat x_\sigma$ from node $\beta$ for
              all $\sigma\in\mathcal{R}^{(\alpha)}_{\text{col},\beta}$;
       \STATE Receive $x|_{\hat\sigma}$ from node $\beta$ for all
              $\sigma\in\mathcal{R}^{(\alpha)}_{\text{col},\beta}$
              that are leaves in $\ctJbeta$
    \ENDFOR;
    \FOR{$\tau\in\ctIalpha$}
      \STATE $\hat y_\tau \gets 0$
    \ENDFOR;
    \FOR{$\beta\in[1:p]$}
      \STATE interaction($\treeroot(\ctalphabeta)$, $\hat x$, $x$,
                         $\hat y$, $y$)
    \ENDFOR;
    \STATE backward($\treeroot(\ctIalpha)$, $\hat y$, $y$)
  \end{algorithmic}
\end{algorithm}

%
% Figure: Runtimes, simple approach
%
\begin{figure}
  \begin{center}
    \includegraphics[width=7cm]{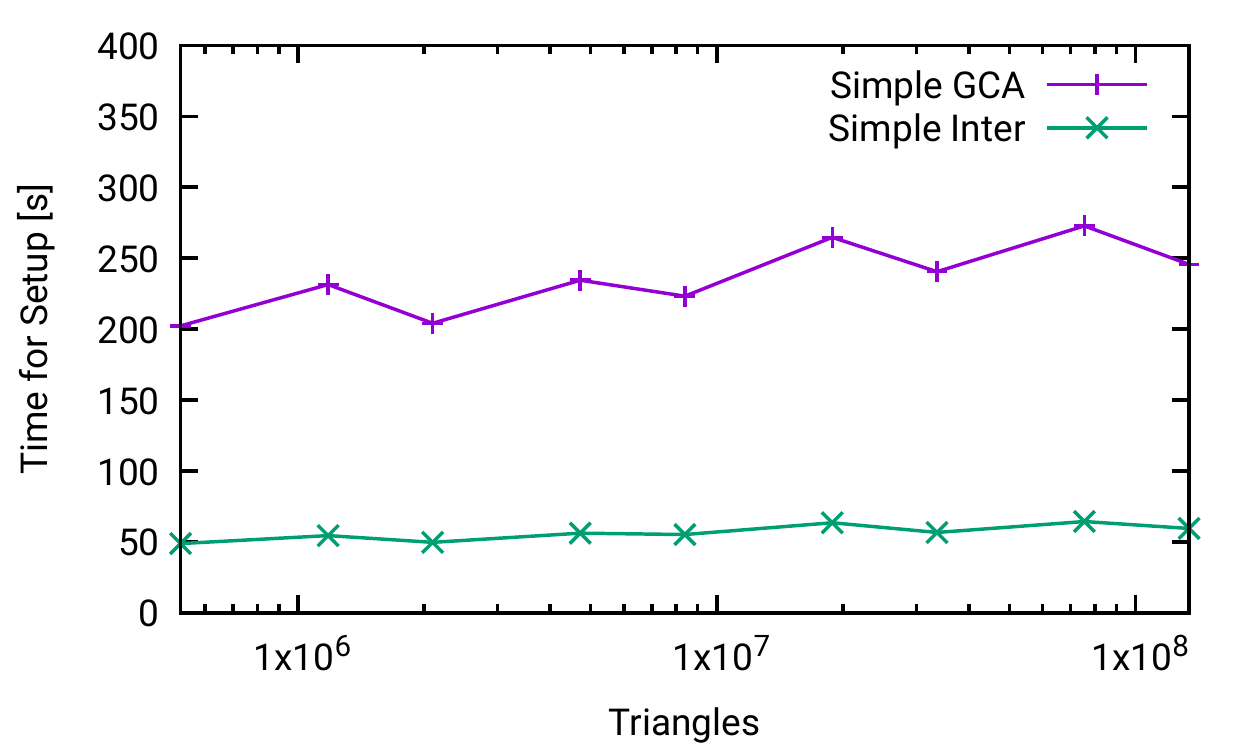}\quad
    \includegraphics[width=7cm]{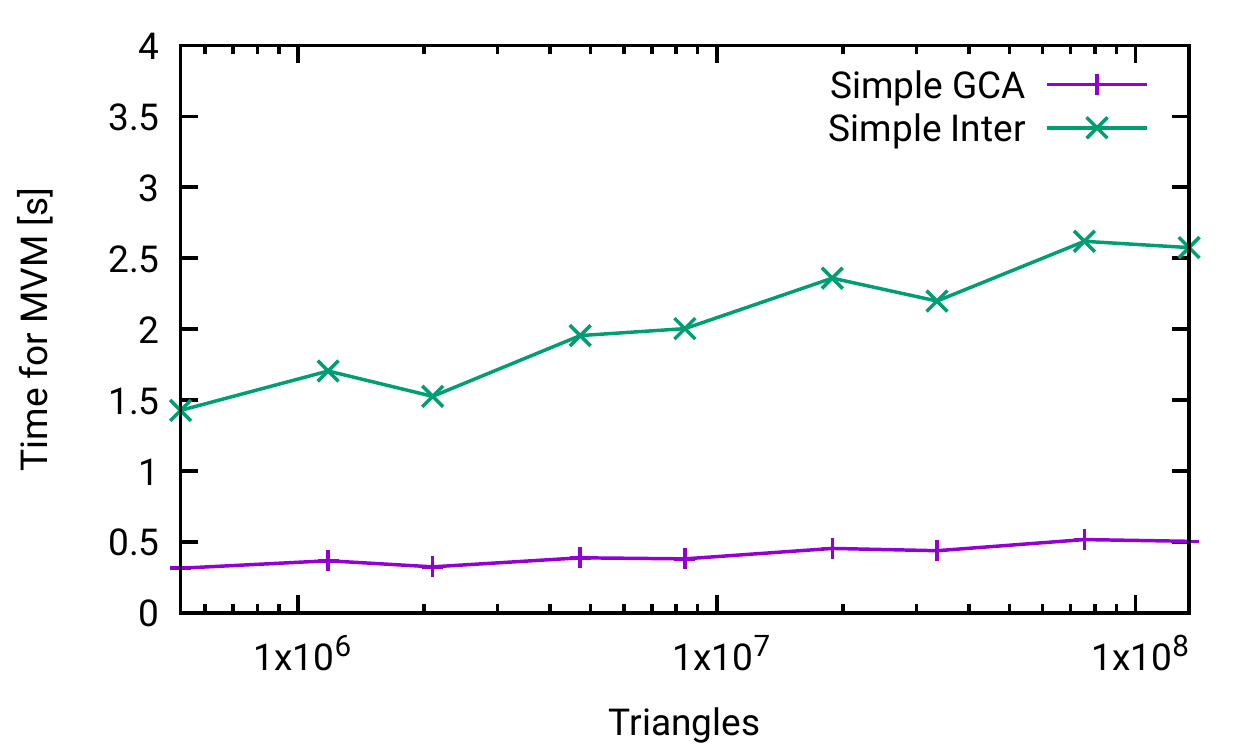}
  \end{center}
  \caption{Runtimes for setup and matrix-vector multiplication for
    $\mathcal{H}^2$-matrices constructed by interpolation and
    GCA-$\mathcal{H}^2$ on HLRN Berlin's ``Lise'' cluster}
  \label{fi:runtimes_simple}
\end{figure}

%
% Table: Runtimes, simple approach
%
\begin{table}
  \caption{Runtimes for setup and matrix-vector multiplication for
    $\mathcal{H}^2$-matrices constructed by interpolation and
    GCA-$\mathcal{H}^2$ on HLRN Berlin's ``Lise'' cluster}
  \label{ta:runtimes_simple}
  \begin{equation*}
    \begin{array}{rr|rr|rr}
        & & \multicolumn{2}{c}{\text{Interpolation}}
          & \multicolumn{2}{|c}{\text{GCA-$\mathcal{H}^2$}}\\
      n & p & \text{Setup} & \text{MVM} & \text{Setup} & \text{MVM}\\
      \hline
      524\,288 & 4 & 48.8 & 1.43 & 202.3 & 0.32 \\
      1\,179\,648 & 8 & 54.5 & 1.70 & 231.4 & 0.37 \\
      2\,097\,152 & 16 & 49.7 & 1.53 & 204.1 & 0.32 \\
      4\,718\,592 & 32 & 56.2 & 1.95 & 234.4 & 0.39 \\
      8\,388\,608 & 64 & 55.2 & 2.00 & 223.2 & 0.38 \\
      18\,874\,368 & 128 & 63.5 & 2.36 & 264.7 & 0.45 \\
      33\,554\,432 & 256 & 56.6 & 2.20 & 240.5 & 0.44 \\
      75\,497\,472 & 512 & 64.4 & 2.62 & 272.7 & 0.52 \\
      134\,217\,728 & 1\,024 & 59.5 & 2.58 & 245.6 & 0.50
    \end{array}
  \end{equation*}
\end{table}

Figure~\ref{fi:runtimes_simple} and Table~\ref{ta:runtimes_simple}
show the runtimes for the matrix setup and the matrix-vector multiplication,
using both interpolation and GCA-$\mathcal{H}^2$ for the approximation.
We approximate the single-layer operator of Laplace's equation
on polyhedral approximations of the unit sphere with piecewise constant basis
functions using Sauter-Schwab quadrature \cite{SA96,SASC11}.
The experiments were run on HLRN Berlin's ``Lise'' cluster,
starting with $4$ MPI processes on one node for $524\,288$
triangles and ending with $1\,024$ MPI processes on $128$ nodes
for $134\,217\,728$ triangles.
For the interpolation approach, we use a constant order of
$m=4$, i.e., cubic tensor polynomials, and an admissibility
parameter of $\eta=1$.
For the GCA-$\mathcal{H}^2$ approach, we also use a constant
quadrature order of $m=4$, but an admissibility parameter
of $\eta=2$ and an accuracy of $10^{-4}$ for the cross approximation.

Both methods appear to achieve an $\mathcal{O}(\log n)$ weak
scaling behaviour both for the setup and the matrix-vector
multiplication, where $n$ is the number of triangles.
Due to $n\sim p$, this is better than expected, since every
node has to set up $\mathcal{O}(n/p + p)$ submatrices.
We can see that interpolation is significantly faster in the
setup phase, while GCA-$\mathcal{H}^2$ is considerably faster
at matrix-vector multiplications.

% ----------------------------------------
% Shared trees
% ----------------------------------------
\section{Shared trees}
\label{se:shared}

By construction, the distributed $\mathcal{H}^2$-matrices introduced
so far have to store at least $p$ blocks, i.e., the blocks
$(\treeroot(\ctIalpha), \treeroot(\ctJbeta))$ for all $\beta\in[1:p]$,
in every node $\alpha\in[1:p]$, therefore the local complexity
can never be below $p k^2$.
This is an unattractive property if we want to treat very large matrices
with very high numbers of nodes.

An elegant approach is to introduce global cluster trees $\ctI$
and $\ctJ$ that have $\ctIalpha$ and $\ctJalpha$, respectively,
as subtrees for all $\alpha\in[1:p]$.
To keep the program's structure close to the original $\mathcal{H}^2$-matrix
construction, it is a good idea to share the additional clusters
among those nodes that own their descendants.
We denote the set of nodes storing a given cluster $\tau\in\ctI$
by $\shareholders(\tau)\subseteq[1:p]$ and require these sets to
have the following properties, similar to the properties of
cluster trees, cf. Definition~\ref{de:cluster_tree}:
\begin{enumerate}
  \item For all $\tau\in\ctIalpha$, $\alpha\in[1:p]$,
    we have $\shareholders(\tau)=\{\alpha\}$.
  \item For all $\tau\in\ctI$ with $\#\shareholders(\tau)>1$ and
    all $\tau_1,\tau_2\in\chil(\tau)$ we have
    \begin{equation*}
      \tau_1\neq\tau_2 \Rightarrow
      \shareholders(\tau_1)\cap\shareholders(\tau_2)=\emptyset.
    \end{equation*}
  \item For all $\tau\in\ctI$ with $\chil(\tau)\neq\emptyset$,
    we have
    \begin{equation*}
      \shareholders(\tau) = \bigcup\{ \shareholders(\tau')
         \ :\ \tau'\in\chil(\tau) \}.
    \end{equation*}
\end{enumerate}
The first condition ensures that the majority of clusters will
only be stored in one node.
The second condition guarantees that of the shared clusters, i.e.,
the clusters with $\#\shareholders(\tau)>1$, every node stores exactly
one per level.
The third condition states that if a node is a shareholder of
a cluster, it is also a shareholder of all of its ancestors.
This property guarantees that every node can find all clusters it
holds a share in by starting from the common root
$\treeroot(\ctI)$.

\begin{figure}
  \begin{center}
    \includegraphics[width=12cm]{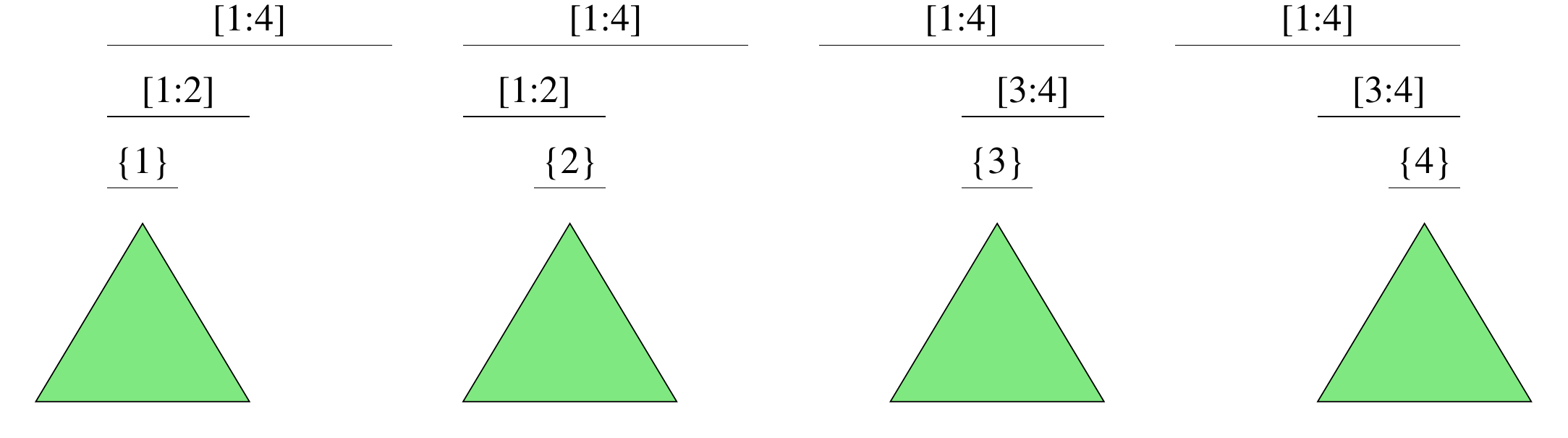}\\[-1cm]
    $\mathcal{T}_{\Idx_1}$\hspace{2.6cm}
    $\mathcal{T}_{\Idx_2}$\hspace{2.6cm}
    $\mathcal{T}_{\Idx_3}$\hspace{2.6cm}
    $\mathcal{T}_{\Idx_4}$\\[0.5cm]
  \end{center}
  \caption{Shared cluster tree: starting with four cluster trees
    $\mathcal{T}_{\Idx_1}$, $\mathcal{T}_{\Idx_2}$, $\mathcal{T}_{\Idx_3}$,
    $\mathcal{T}_{\Idx_4}$, we add a new shared cluster $\tau_{1,2}$ with
    $\chil(\tau_{1,2})=\{\treeroot(\mathcal{T}_{\Idx_1}),
    \treeroot(\mathcal{T}_{\Idx_2})\}$ and $\shareholders(\tau_{1,2})=[1:2]$,
    a new shared cluster $\tau_{3,4}$ with
    $\chil(\tau_{3,4})=\{\treeroot(\mathcal{T}_{\Idx_3}),
    \treeroot(\mathcal{T}_{\Idx_4})\}$ and $\shareholders(\tau_{3,4})=[3:4]$,
    and a new shared root $\tau_{1,4}$ with
    $\chil(\tau_{1,4})=\{\tau_{1,2},\tau_{3,4}\}$
    and $\shareholders(\tau_{3,4})=[1:4]$}
  \label{fi:shared_cluster}
\end{figure}

One approach to constructing a shared cluster tree is to
choose a characteristic point $x_\alpha\in\Gamma_\alpha$ for
every $\alpha\in[1:p]$ and then apply standard clustering techniques,
e.g., geometric bisection, to construct a cluster tree for these
points.
The shareholders for one of the clusters are then simply the
nodes $\alpha\in[1:p]$ corresponding to the characteristic points
contained in a cluster.
Once a cluster contains only one point $x_\alpha$, we can attach
the tree $\ctIalpha$.
An example for $p=4$ nodes can be found in Figure~\ref{fi:shared_cluster}.

%
% Figure: Shared cluster tree
%

We can use the same approach to construct a shared cluster tree
$\ctJ$ from $\ctJalpha$ for $\alpha\in[1:p]$, and starting with
the shared trees $\ctI$ and $\ctJ$, we can obtain a shared
block tree $\ctIJ$ by using Algorithm~\ref{al:build_block} with
a few modifications.
The most important change is that we make communication more
efficient by appointing a ``manager'' to each cluster, i.e.,
we choose
\begin{align*}
  \manager(\tau) &\in \shareholders(\tau) &
  &\text{ for all } \tau\in\ctI,\\
  \manager(\sigma) &\in \shareholders(\sigma) &
  &\text{ for all } \sigma\in\ctJ,
\end{align*}
and follow the convention that only the manager of a cluster is
allowed to send and receive information regarding this cluster.
If the other shareholders need the information received by the
manager, it can be provided efficiently by a broadcast operation.

The second change follows from the fact that in the shared
cluster tree, the children of a cluster may have different
shareholders than the parent.
Assume that we have a block $(\tau,\sigma)$ that is not
admissible with $\#\chil(\tau)>1$ and $\#\chil(\sigma)>1$.

If $\#\shareholders(\tau)>1$ holds, the manager $\beta:=\manager(\sigma)$
of $\sigma$ has to send information on the children of $\sigma$ to the
\emph{different} managers $\alpha':=\manager(\tau')$ of the children
$\tau'\in\chil(\tau)$ of $\tau$.

If $\#\shareholders(\sigma)>1$ holds, the manager $\alpha:=\manager(\tau)$
of $\tau$ has to send information on the children of $\tau$ to the
\emph{different} managers $\beta':=\manager(\sigma')$ of the children
$\sigma'\in\chil(\sigma)$ of $\sigma$.

We can solve this issue by modifying the structure of the send
trees:
shared send trees are no longer subtrees of the corresponding
cluster trees, but consist of pairs of a cluster and the
manager receiving this cluster.

%
% Definition: Shared send tree
%
\begin{definition}[Shared send tree]
Let $\ctI$ and $\ctJ$ be shared cluster trees, and let $\ctIJ$
be the resulting block tree.
The corresponding \emph{shared send trees} $\mathcal{S}^{(\alpha)}_\text{row}$
and $\mathcal{S}^{(\beta)}_\text{col}$, $\alpha,\beta\in[1:p]$, are defined
as follows:
\begin{itemize}
  \item For every $\alpha\in[1:p]$, every vertex
    $t\in\mathcal{S}^{(\alpha)}_\text{row}$ is of
    the form $t=(\tau,\beta)$ with a row cluster $\tau\in\ctI$ and
    a column manager $\beta\in[1:p]$.

    For every $\beta\in[1:p]$, every vertex
    $s\in\mathcal{S}^{(\beta)}_\text{col}$ is of
    the form $s=(\sigma,\alpha)$ with a column cluster $\sigma\in\ctJ$ and
    a row manager $\alpha\in[1:p]$.
  \item For every $\alpha\in[1:p]$, the root of
    $\mathcal{S}^{(\alpha)}_\text{row}$ is
    $(\treeroot(\ctI),\manager(\treeroot(\ctJ)))$.
  
    For every $\beta\in[1:p]$, the root of
    $\mathcal{S}^{(\alpha)}_\text{col}$ is
    $(\treeroot(\ctJ),\manager(\treeroot(\ctI)))$.
  \item A vertex $t=(\tau,\beta)\in\mathcal{S}^{(\alpha)}_\text{row}$,
    $\alpha\in[1:p]$, has children if and only if $\chil(\tau)\neq\emptyset$
    and there is a $\sigma\in\ctJ$ such that $(\tau,\sigma)\in\ctIJ$ is
    inadmissible and $\beta=\manager(\sigma)$.
    Then the children are
    \begin{align*}
      \chil(t,\beta) &=
                   \{ (\tau',\beta')\ :\ \tau'\in\chil(\tau),
                       \ \beta'=\manager(\sigma'),
                       \ \sigma'\in\chil(\sigma) \}
    \intertext{if $\chil(\sigma)\neq\emptyset$ and otherwise}
      \chil(t,\beta) &=
                   \{ (\tau',\beta)\ :\ \tau'\in\chil(\tau) \}.
    \end{align*}
  \item A vertex $s=(\sigma,\alpha)\in\mathcal{S}^{(\beta)}_\text{row}$,
    $\beta\in[1:p]$, has children if and only if $\chil(\sigma)\neq\emptyset$
    and there is a $\tau\in\ctJ$ such that $(\tau,\sigma)\in\ctIJ$ is
    inadmissible and $\alpha=\manager(\tau)$.
    Then the children are
    \begin{align*}
      \chil(s,\alpha) &= 
                   \{ (\sigma',\alpha')\ :\ \sigma'\in\chil(\sigma),
                       \ \alpha'=\manager(\tau'),
                       \ \tau'\in\chil(\tau) \}
    \intertext{if $\chil(\tau)\neq\emptyset$ and otherwise}
      \chil(s,\alpha) &=
                   \{ (\sigma',\alpha)\ :\ \sigma'\in\chil(\sigma) \}.
    \end{align*}
\end{itemize}
\end{definition}

In short, an entry $s=(\sigma,\alpha)\in\mathcal{S}^{(\beta)}_\text{col}$
means that the node $\beta$ has to send information on the cluster $\sigma$
to the node $\alpha$, while an entry
$t=(\tau,\beta)\in\mathcal{S}^{(\alpha)}_\text{row}$ means that the
node $\alpha$ has to send information on the cluster $\tau$ to
the node $\beta$.

Algorithm~\ref{al:build_block} can be adjusted to handle shared
cluster trees, this results in Algorithm~\ref{al:build_shared},
where some intermediate steps have been abbreviated:
``Expect to receive'' and ``expect to send'' means adding the
corresponding clusters to suitable sets similar to
$\mathcal{G}_{\text{row},\beta}$ and $\mathcal{P}_{\text{row},\beta}$
in Algorithm~\ref{al:build_block}, with the purpose of concentrating
communication operations at the ``Send and receive expected children''
point, where they can be carried out, e.g., by MPI collective communication
functions.

%
% Algorithm: Construction of block, send, and receive trees
%
\begin{algorithm}
  \caption{Construction of shared block trees, send and receive trees}
  \label{al:build_shared}
  \begin{algorithmic}[1]
    \STATE \textbf{procedure} build\_shared(
               $\ctI$, $\ctJ$, \textbf{var} $\ctIJ$,
               $\mathcal{R}^{(\alpha)}_\text{row}$,
               $\mathcal{R}^{(\alpha)}_\text{col}$,
               $\mathcal{S}^{(\alpha)}_\text{row}$,
               $\mathcal{S}^{(\alpha)}_\text{col}$);
    \STATE $\treeroot(\ctIJ)
             \gets (\treeroot(\ctI),\treeroot(\ctJ))$;\ 
           $\mathcal{A}_\text{row} \gets \{ \treeroot(\ctIJ) \}$;\ 
           $\mathcal{A}_\text{col} \gets \{ \treeroot(\ctIJ) \}$;
    \STATE $\treeroot(\mathcal{R}^{(\alpha)}_{\text{row}})
             \gets \treeroot(\ctI)$;\ 
           $\treeroot(\mathcal{R}^{(\alpha)}_{\text{col}})
             \gets \treeroot(\ctJ)$;
    \STATE $\treeroot(\mathcal{S}^{(\alpha)}_{\text{row}})
             \gets (\treeroot(\ctI), \manager(\treeroot(\ctJ)))$;
    \STATE $\treeroot(\mathcal{S}^{(\alpha)}_{\text{col}})
             \gets (\treeroot(\ctJ), \manager(\treeroot(\ctI)))$;
    \REPEAT
      \FOR{$(\tau,\sigma)\in\mathcal{A}_\text{row}$}
        \IF{$(\tau,\sigma)$ not admissible, $\chil(\sigma)\neq\emptyset$
            and $\alpha = \manager(\tau)$}
          \STATE Expect to receive the children of $\sigma$
                 from $\manager(\sigma)$;
          \STATE Expect to send the children of $\tau$ to
                 $\manager(\sigma)$
        \ENDIF
      \ENDFOR
      \FOR{$(\tau,\sigma)\in\mathcal{A}_\text{col}$}
        \IF{$(\tau,\sigma)$ not admissible, $\chil(\tau)\neq\emptyset$
            and $\alpha = \manager(\sigma)$}
          \STATE Expect to receive the children of $\tau$
                 from $\manager(\tau)$;
          \STATE Expect to send the children of $\sigma$ to
                 $\manager(\tau)$
        \ENDIF
      \ENDFOR;
      \STATE Send and receive expected children;
      \STATE $\mathcal{A}_\text{row}^\text{old} \gets \mathcal{A}_\text{row}$;
             \ $\mathcal{A}_\text{row} \gets \emptyset$;
      \FOR{$(\tau,\sigma)\in\mathcal{A}_\text{row}^\text{old}$}
        \IF{$(\tau,\sigma)$ not admissible, $\chil(\sigma)\neq\emptyset$}
          \STATE $\manager(\tau)$ sends the children of $\sigma$
                 to all nodes in $\shareholders(\tau)$;
          \STATE Create the children of the block $(\tau,\sigma)$;
          \STATE Add $(\tau',\manager(\sigma'))$ to
                 $\chil(\tau,\manager(\sigma))$ in
                 $\mathcal{S}^{(\alpha)}_\text{row}$ for
                 $(\tau',\sigma')\in\chil(\tau,\sigma)$;
          \STATE Add $\sigma'$ to $\chil(\sigma)$ in
                 $\mathcal{R}^{(\alpha)}_\text{col}$ for
                 $(\tau',\sigma')\in\chil(\tau,\sigma)$;
          \STATE Add $(\tau',\sigma')$ to $\mathcal{A}_\text{row}$
                 for all $(\tau',\sigma')\in\chil(\tau,\sigma)$
                 with $\alpha\in\shareholders(\tau')$
        \ENDIF
      \ENDFOR
      \STATE $\mathcal{A}_\text{col}^\text{old} \gets \mathcal{A}_\text{col}$;
             \ $\mathcal{A}_\text{col} \gets \emptyset$;
      \FOR{$(\tau,\sigma)\in\mathcal{A}_\text{col}^\text{old}$}
        \IF{$(\tau,\sigma)$ not admissible, $\chil(\tau)\neq\emptyset$}
          \STATE $\manager(\sigma)$ sends the children of $\tau$
                 to all nodes in $\shareholders(\sigma)$;
          \STATE Create the children of the block $(\tau,\sigma)$;
          \STATE Add $(\sigma',\manager(\tau'))$ to
                 $\chil(\sigma,\manager(\tau))$ in
                 $\mathcal{S}^{(\alpha)}_\text{col}$ for
                 $(\tau',\sigma')\in\chil(\tau,\sigma)$;
          \STATE Add $\tau'$ to $\chil(\tau)$ in
                 $\mathcal{R}^{(\alpha)}_\text{row}$ for
                 $(\tau',\sigma')\in\chil(\tau,\sigma)$;
          \STATE Add $(\tau',\sigma')$ to $\mathcal{A}_\text{col}$
                 for all $(\tau',\sigma')\in\chil(\tau,\sigma)$
                 with $\alpha\in\shareholders(\sigma')$
        \ENDIF
      \ENDFOR
      \STATE busy=$(\mathcal{A}_\text{row}\neq\emptyset \text{ or }
                    \mathcal{A}_\text{col}\neq\emptyset)$;\ 
             reduce\_or\_bool(busy)
    \UNTIL{busy=false}
  \end{algorithmic}
\end{algorithm}

Algorithm~\ref{al:mvm} for the matrix-vector multiplication also has
to be adjusted in order to handle shared cluster trees:
in the forward and backward transformation, only the manager of a
cluster is allowed to perform arithmetic operations, therefore
if one of the children of a cluster has a different manager than
the parent, an additional communication operation is required.
The resulting shared forward and backward transformations are
given in Algorithms~\ref{al:forward_shared} and \ref{al:backward_shared}.
In the interaction phase, only the manager of the row cluster of
a block is allowed to perform arithmetic operations.
Algorithm~\ref{al:mvm_shared} uses the shared send and receive
trees to ensure that the coefficients $\hat x_\sigma$ are properly
transmitted from the node managing $\sigma$ to all nodes that need
these coefficients.

%
% Algorithm: Forward transformation
%
\begin{algorithm}
  \caption{Shared forward transformation}
  \label{al:forward_shared}
  \begin{algorithmic}[1]
    \STATE \textbf{procedure} forward\_shared($\sigma$, $x$,
               \textbf{var} $\hat x$);
    \IF{$\chil(\sigma)=\emptyset$}
      \STATE $\hat x_\sigma \gets W_\sigma^T x|_{\hat\sigma}$
    \ELSE
      \STATE $\hat x_\sigma \gets 0$;
      \FOR{$\sigma'\in\chil(\sigma)$}
        \IF{$\alpha\in\shareholders(\sigma')$}
          \STATE forward($\sigma'$, $x$, $\hat x$)
        \ENDIF;
        \IF{$\alpha=\manager(\sigma)$ and $\alpha\neq\manager(\sigma')$}
          \STATE Receive $\hat x_{\sigma'}$ from node $\alpha'=\manager(\sigma')$
        \ELSIF{$\alpha\neq\manager(\sigma)$ and $\alpha=\manager(\sigma')$}
          \STATE Send $\hat x_{\sigma'}$ to node $\alpha'=\manager(\sigma)$
        \ENDIF;
        \STATE $\hat x_\sigma \gets \hat x_\sigma + F_{\sigma'}^T \hat x_{\sigma'}$
      \ENDFOR
    \ENDIF
  \end{algorithmic}
\end{algorithm}

%
% Algorithm: Shared backward transformation
%
\begin{algorithm}
  \caption{Shared backward transformation}
  \label{al:backward_shared}
  \begin{algorithmic}[1]
    \STATE \textbf{procedure} backward\_shared($\tau$,
               \textbf{var} $\hat y$, $y$);
    \IF{$\chil(\tau)=\emptyset$}
      \STATE $y|_{\hat\tau}\gets y|_{\hat\tau} + V_\tau \hat y_\tau$
    \ELSE
      \FOR{$\tau'\in\chil(\tau)$}
        \STATE $\hat y_{\tau'} \gets \hat y_{\tau'} + E_{\tau'} \hat y_\tau$;
        \IF{$\alpha=\manager(\tau)$ and $\alpha\neq\manager(\tau')$}
          \STATE Send $\hat y_{\tau'}$ to node $\alpha'=\manager(\tau')$
        \ELSIF{$\alpha\neq\manager(\tau)$ and $\alpha=\manager(\tau')$}
          \STATE Receive $\tilde y_{\tau'}$ from node $\alpha'=\manager(\tau)$
                 and add to $\hat y_{\tau'}$
        \ENDIF;
        \IF{$\alpha\in\shareholders(\tau')$}
          \STATE backward($\tau'$, $\hat y$, $y$)
        \ENDIF
      \ENDFOR
    \ENDIF
  \end{algorithmic}
\end{algorithm}

%
% Algorithm: Shared matrix-vector multiplication
%
\begin{algorithm}
  \caption{Shared matrix-vector multiplication}
  \label{al:mvm_shared}
  \begin{algorithmic}[1]
    \STATE \textbf{procedure} mvm\_shared(
               $\mathcal{T}_{\Idx\times\Jdx}$,
               $\mathcal{R}^{(\alpha)}_\text{col}$,
               $\mathcal{S}^{(\alpha)}_\text{col}$,
               $x$, \textbf{var} y);
    \STATE forward\_shared($\treeroot(\ctJ)$, $x$, $\hat x$);
    \FOR{$(\sigma,\beta)\in\mathcal{S}^{(\alpha)}_\text{col}$}
       \STATE Send $\hat x_\sigma$ to node $\beta$;
       \STATE Send $x|_{\hat\sigma}$ to node $\beta$
              if $\chil(\sigma)=\emptyset$
    \ENDFOR;
    \FOR{$\sigma\in\mathcal{R}^{(\alpha)}_\text{col}$}
       \STATE Receive $\hat x_\sigma$ from node $\beta=\manager(\sigma)$;
       \STATE Receive $x|_{\hat\sigma}$ from node $\beta=\manager(\sigma)$
              if $\chil(\sigma)=\emptyset$
    \ENDFOR;
    \FOR{$\tau\in\ctIalpha$}
      \STATE $\hat y_\tau \gets 0$
    \ENDFOR;
    \STATE interaction\_shared($\treeroot(\ctIJ)$, $\hat x$, $x$,
                       $\hat y$, $y$);
    \STATE backward\_shared($\treeroot(\ctI)$, $\hat y$, $y$)
  \end{algorithmic}
\end{algorithm}

%
% Figure: Runtimes, simple and shared approach
%
\begin{figure}
  \begin{center}
    \includegraphics[width=7cm]{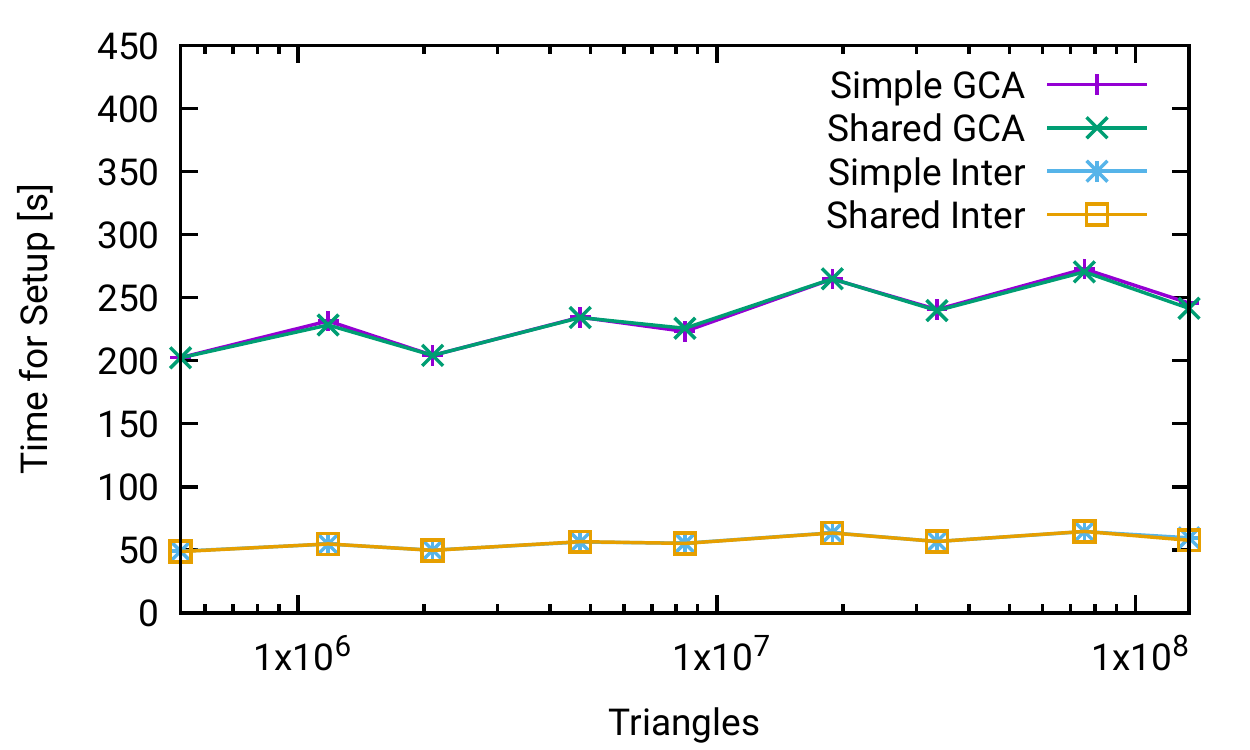}\quad
    \includegraphics[width=7cm]{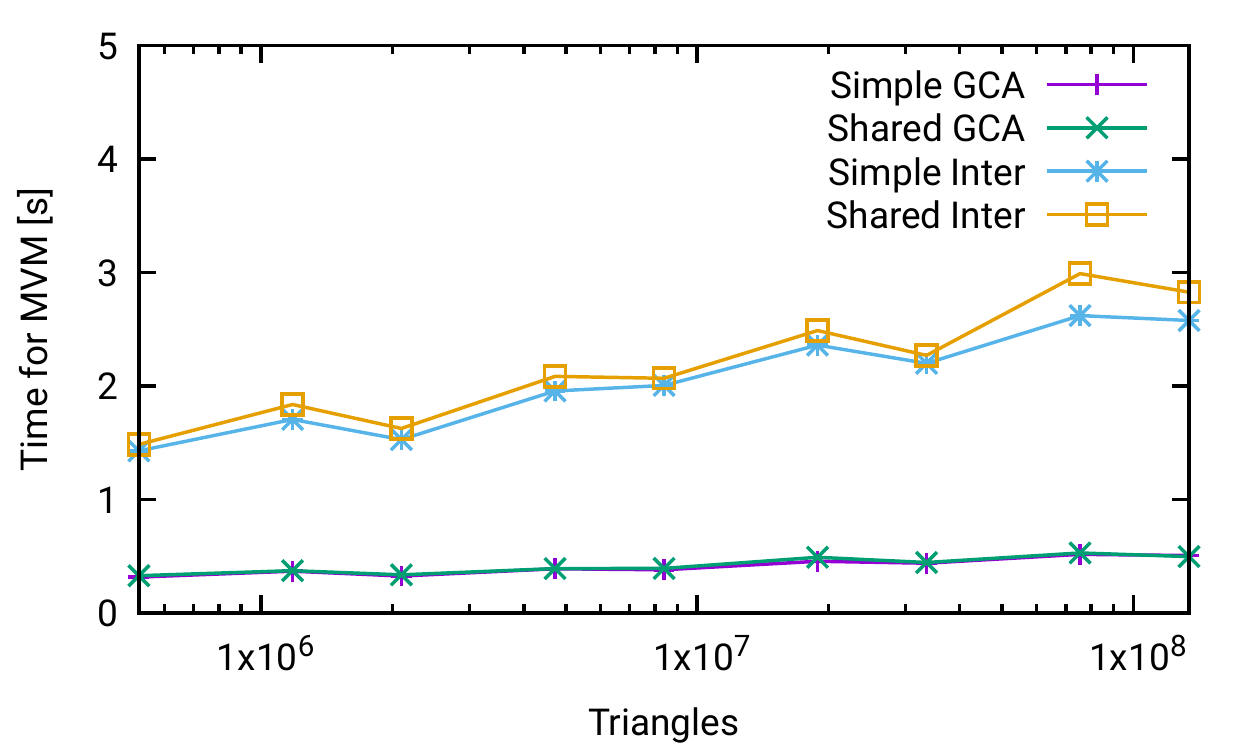}
  \end{center}
  \caption{Runtimes for simple and shared setup and matrix-vector
    multiplication for $\mathcal{H}^2$-matrices constructed by
    interpolation and GCA-$\mathcal{H}^2$ on HLRN Berlin's ``Lise'' cluster}
  \label{fi:runtimes_shared}
\end{figure}

%
% Table: Runtimes, simple and shared approach
%
\begin{table}
  \caption{Runtimes for setup and matrix-vector multiplication for
    $\mathcal{H}^2$-matrices constructed by interpolation and
    GCA-$\mathcal{H}^2$ on HLRN Berlin's ``Lise'' cluster}
  \label{ta:runtimes_shared}
  \begin{equation*}
    \begin{array}{rr|rr|rr}
        & & \multicolumn{2}{c}{\text{Interpolation}}
          & \multicolumn{2}{|c}{\text{GCA-$\mathcal{H}^2$}}\\
      n & p & \text{Setup} & \text{MVM} & \text{Setup} & \text{MVM}\\
      \hline
      524\,288 & 4 & 48.6 & 1.48 & 202.2 & 0.33 \\
      1\,179\,648 & 8 & 54.6 & 1.84 & 228.3 & 0.37 \\
      2\,097\,152 & 16 & 49.5 & 1.62 & 204.2 & 0.34 \\
      4\,718\,592 & 32 & 56.5 & 2.08 & 234.1 & 0.39 \\
      8\,388\,608 & 64 & 55.0 & 2.07 & 225.6 & 0.39 \\
      18\,874\,368 & 128 & 63.3 & 2.49 & 264.8 & 0.49 \\
      33\,554\,432 & 256 & 56.6 & 2.27 & 239.6 & 0.44 \\
      75\,497\,472 & 512 & 64.5 & 2.99 & 270.3 & 0.53 \\
      134\,217\,728 & 1\,024 & 57.6 & 2.83 & 241.4 & 0.50
    \end{array}
  \end{equation*}
\end{table}

Figure~\ref{fi:runtimes_shared} and Table~\ref{ta:runtimes_shared}
show that the shared approach does not lead to significant changes
in the runtimes for the setup phase, while it is even a little
slower for the matrix-vector multiplication in case interpolation
is used.
This may be due to the implementation: currently MPI collective
communication operations are used to transmit the coefficients
obtained in the forward transformation, and it is to be expected
that these operations take at least $\Omega(\log p)$ operations
per level, even if the shared algorithms transfer data only between
a small number of nodes.
Answering the question whether using simple point-to-point
communication operations leads to a significant improvement will
be the topic future research.

Another topic will be the inclusion of more advanced algorithms,
e.g., the different phases of the on-the-fly recompression algorithm
\cite[Section~6.6]{BO10} will require us to compute different weight
matrices following the patterns of the forward transformation and
the interaction phase.

% ----------------------------------------
% Bibliography
% ----------------------------------------
\bibliographystyle{plain}
\bibliography{hmatrix}

\end{document}